\documentclass[amscd,amssymb,epic,eepic,leqno,12pt]{amsart}
\usepackage{latexsym,amsmath,amsfonts,amscd,amssymb}
\newenvironment{pf}{\it Proof:\rm}{\qed}

\newcommand{\AAA}{{\mathcal A}}
\newcommand{\CCC}{{\mathcal C}}
\newcommand{\EEE}{{\mathcal E}}
\newcommand{\FFF}{{\mathcal F}}
\newcommand{\GGG}{{\mathcal G}}
\newcommand{\GGGS}{{\mathcal G}_{sh}}
\newcommand{\HHHS}{{\mathcal H}_{sh}}
\newcommand{\GGC}{{\mathcal G^{\CC}}}
\newcommand{\HHH}{{\mathcal H}}
\newcommand{\HHC}{{\mathcal H^{\CC}}}

\newcommand{\OOO}{{\mathcal O}}
\newcommand{\Ok}{{\mathcal O^{k}}}

\newcommand{\TTT}{{\mathcal T}}

\newcommand{\VVV}{{\mathcal V}}

\newcommand{\XXX}{{\mathcal X}}

\newcommand{\Met}{{\mathcal M}\!et^p_2}
\newcommand{\MetB}{{\mathcal M}\!et^p_{2,B}}

\newcommand{\noi}{\noindent}
\newcommand{\la}{\langle}
\newcommand{\ra}{\rangle}

\newtheorem{theorem}{\bf Theorem}[section]
\newtheorem{assumption}[theorem]{\bf Assumptions}
\newtheorem{lemma}[theorem]{\bf Lemma}
\newtheorem{prop}[theorem]{\bf Proposition}
\newtheorem{corollary}[theorem]{\bf Corollary}
\newtheorem{definition}[theorem]{\bf Definition}
\newtheorem{example}[theorem]{\bf Example}

\theoremstyle{remark}
\newtheorem{remark}[theorem]{\bf Remark}
\newcommand{\be}{\begin{equation}}
\newcommand{\ee}{\end{equation} }

\newfont{\bfc}{cmbsy10 scaled 1200}  
\newfont{\dr}{msbm10 scaled \magstep1}  
\newfont{\sdr}{msbm8}  
\newfont{\gl}{eufm10 scaled \magstep1}  
\DeclareFontFamily{OT1}{rsfs}{}
\DeclareFontShape{OT1}{rsfs}{n}{it}{<->rsfs10}{}
\newcommand{\V}{{\mathbb V}}
\newcommand{\VV}{{V}}
\newcommand{\W}{{W}}
\newcommand{\Y}{{\mathbb Y}}
\newcommand{\YY}{{Y}}
\newcommand{\CC}{{\mathbb C}}

\newcommand{\RR}{{\mathbb R}}
\newcommand{\WW}{{\mathbb W}}
\newcommand{\ZZ}{{\mathbb Z}}

\newcommand{\glie}{{\mathfrak g}}
\newcommand{\hlie}{{\mathfrak h}}
\newcommand{\klie}{{\mathfrak k}}

\newcommand{\ulie}{{\mathfrak u}}

\newcommand{\Ad}{\operatorname{Ad}}
\newcommand{\ad}{\mathop{{\fam0 ad}}\nolimits}

\newcommand{\codim}{\operatorname{codim}}

\newcommand{\End}{\operatorname{End}}
\newcommand{\GL}{\operatorname{GL}}

\newcommand{\Hom}{\operatorname{Hom}}
\newcommand{\Id}{\operatorname{Id}}

\newcommand{\Ker}{\operatorname{Ker}}
\newcommand{\Lie}{\operatorname{Lie}}

\newcommand{\rk}{\operatorname{rk}}

\newcommand{\Tr}{\operatorname{Tr}}
\newcommand{\U}{\operatorname{U}}

\newcommand{\Vol}{\operatorname{Vol}}
\renewcommand{\exp}{\operatorname{exp}}

\newcommand{\ov}{\overline}

\newcommand{\imag}{\sqrt{-1}}
\newcommand{\qu}{/\kern-.7ex/}
\newcommand{\exh}{\to\kern-1.8ex\to}
\newcommand{\wo}{\widehat{\otimes}}
\newcommand{\rank}{{\mathrm {rank}}}

\begin{document}
\begin{titlepage}
\noindent
{\Large\textbf{Relative Hitchin--Kobayashi correspondences\\ for
principal pairs}}
\bigskip

\noindent
\textbf{Steven B. Bradlow}\footnotemark[1]$^{,}$\footnotemark[2] \\
Department of Mathematics, \\ University of Illinois, \\ Urbana, IL
61801, USA \\ {\it E-mail:}  \texttt{bradlow@math.uiuc.edu}
\medskip

\noindent
\textbf{Oscar Garc{\'\i}a-Prada}\footnotemark[1]$^{,}$\footnotemark[3]$^{,}$\footnotemark[4] \\
Departamento de Matem{\'a}ticas, \\ Universidad Aut{\'o}noma de
Madrid, \\ 28049 Madrid, Spain \\ {\it E-mail:}
\texttt{oscar.garcia-prada@uam.es}

\medskip
\noindent
\textbf{Ignasi Mundet i Riera}\footnotemark[1]$^{,}$\footnotemark[3]$^{,}$\footnotemark[5] \\
Departamento de Matem{\'a}ticas, \\ Universidad Aut{\'o}noma de
Madrid, \\ 28049 Madrid, Spain \\ {\it Current Address:}\\ Dep.
Matem\`atica Aplicada I ETSEIB, UPC \\ c/Diagonal 647 08028
Barcelona, Spain \\ {\it E-mail:}
\texttt{ignasi.mundet@uam.es}
\date{June  1, 2002}
\footnotetext[1]{
Members of VBAC (Vector Bundles on Algebraic Curves), which is
partially supported by EAGER (EC FP5 Contract no.\
HPRN-CT-2000-00099) and by EDGE (EC FP5 Contract no.\
HPRN-CT-2000-00101).}
\footnotetext[2]{Partially supported by
the National Science Foundation under grant DMS-0072073. }

\footnotetext[3]{Partially supported by the Ministerio de
Ciencia y Tecnolog\'{\i}a (Spain) under grant BFM2000-0024.}

\footnotetext[4]{Partially supported by a British EPSRC grant
(October-December 2001).}

\footnotetext[5]{Partially supported by an
EDGE  grant (October-December 2001). }

\bigskip

\vfill

\paragraph{\bf Abstract.} A principal pair consists of a holomorphic principal
$G$-bundle together with a holomorphic section of an associated Kaehler
fibration. Such objects support natural gauge theoretic equations
coming from a moment map condition, and also admit a notion of
stability based on Geometric Invariant Theory. The Hitchin--Kobayashi
correspondence for principal pairs identifies stability as the
condition for the existence of solutions to the equations. In this
paper we generalize these features in a way which allows the full
gauge group of the principal bundle to be replaced by certain proper
subgroups. Such a generalization is needed in order to use principal
pairs as a general framework for describing augmented holomorphic
bundles. We illustrate our results with applications to well known
examples.

\vfill
\newpage
\end{titlepage}

\section{Introduction}
The Hitchin--Kobayashi correspondence relates a notion of stability
inspired by Geometric Invariant Theory and a set of partial
differential equations coming from Gauge Theory.  Originally proven
for holomorphic bundles (cf. \cite{NS,D,UY,RS}), where it relates slope
stability to the Hermitian-Einstein equations,  similar
correspondences are known to hold in more general settings. Despite
the uniformity of the pattern, each instance of the correspondence
has required individualized treatment. Such ad hoc approaches
inevitably obscure the underlying general principles. There is thus a
clear need for a unifying framework which clarifies the origins of
the common features.

The full scope of the correspondences is determined, at least in the
current state of the art, by the collection of instances where they
have been established. In all cases, the setting consists of
holomorphic bundles with some extra prescribed holomorphic data.
Called augmented holomorphic bundles, a useful illustrative example
is that of the holomorphic pair $(E,\phi)$ in which $E$ is a
holomorphic bundle and $\phi\in H^0(E)$\ is a holomorphic section.
Other examples which have emerged naturally in other settings or have
turned out to have interesting applications are surveyed in
\cite{BDGW} and also in \cite{Sch,M,AG}.

A general framework intended to unify a broad class of augmented
bundles was introduced by one of us in \cite{M}, building on ideas
introduced in \cite{Ba}. The key elements in this framework, called
principal pairs,  are objects defined by (i) a principal $K$-bundle
$P_K\longrightarrow X$, where
$X$ is a compact Kaehler manifold and $K$ is the compact real form of
a complex Lie group $G$, and (ii) a Kaehler manifold, $\Y$,\ with a
Hamiltonian
$K$-action. A principal pair is then a pair $(A,\Phi)$, where $A$\ is a
connection on $P_K$, and $\Phi$ is a section of the associated fibre
bundle $\YY_K=P_K\times_K\Y$. Equivalently, one can replace
$P_K$ by the principal $G$-bundle $P_G=P_K\times_KG$ and consider the
unique extension of the $K$-action to a holomorphic $G$-action on $\Y$
(this exists because the complex structure on $\Y$ is integrable).
The principal pair is then described by a holomorphic structure on
$P_G$ together with a section of the associated bundle
$\YY_G=P_G\times_G\Y$.

As shown in [M], there is a natural Hitchin--Kobayashi correspondence
for such pairs. Moreover, many familiar instances of
Hitchin--Kobayashi correspondences on augmented bundles follow as
special cases of the general result for principal pairs. There are
however important examples, including coherent systems, and Higgs
bundles, which do not fit this mold.

In each of these cases where the principal pairs framework turns out
to be inadequate, the reason is the same. By an appropriate choice of
principal bundle and Kaehler manifold, the augmented bundles can
indeed be described as particular types of principal pairs. The
problem is that the set of all augmented bundles of the given type
corresponds to only a subset of all the corresponding principal
pairs. Furthermore, the automorphism groups for the augmented bundles
correspond to only a subgroup of the principal pairs automorphism
groups.

Our goal in this paper is to overcome these shortcomings in the
principal pairs framework.  Since the automorphism group for a
principal pair corresponds to the gauge group of the principal
bundle, we are led to reexamine the role played by the gauge group in
the theory of principal pairs. In our main result (Theorem
\ref{theorem:main} in Section \ref{sect:main})
we show that subject to certain natural constraints, the results of
[M] can indeed be generalized to allow for the full gauge group to
replaced by a subgroup.

As in the proof given in [M], our proof is motivated by the analogy
between the  Hitchin--Kobayashi correspondence and the results of
Kirwan and Kempf--Ness (see Remark \ref{rk:restriction} for an
important subtlety which explains why in some cases this is only an
approximate analogy). The Kirwan--Kempf--Ness results apply to finite
dimensional (Geometric Invariant Theory or Kaehler) quotients, where
they relate an appropriate notion of stability to the vanishing of a
symplectic moment map. The key to understanding this relation comes
from the re-formulation of stability as a numerical condition
expressed in terms of a function called the Hilbert numerical
function. This function, in turn, can be expressed in terms of an
integral of the moment map.

In the setting which describes principal pairs, the required quotient
construction involves the action of an infinite dimensional (gauge)
group on the infinite dimensional configuration space of all
holomorphic pairs. One can no longer directly apply Geometric
Invariant Theory, but, as shown in [M], one can still construct an
integral of the moment map. One can thus define stability in
precisely the same way as dictated by the Hilbert numerical
condition.  The Hitchin--Kobayashi correspondence between stable
orbits of the complex gauge group and orbits which contain a zero of
the moment map is then revealed as a precise analog of the finite
dimensional correspondence.

The generalization we require is easily accomplished in the {\it
finite} dimensional setting, where it corresponds to replacing the
action of a reductive group on a projective variety by the action of
a reductive subgroup. In that setting it is clear how to modify the
GIT notion of stability, and also how to modify the appropriate
moment map. Our main task is thus to show how to make the analogous
modifications in the {\it infinite} dimensional gauge theoretic
setting. As described in Section \ref{sect:main}, the resulting
generalized Hitchin--Kobayashi correspondence still relates zeros of a
moment map and a notion of stability, but now both are defined with
respect to the action of a subgroup of a gauge group.

In order to check that the resulting equations and notions of
stability correspond to those obtained by the ad hoc methods used on
specific examples, it is useful to understand how they relate to
their counterparts for the full gauge group. For both the equations
and the stability conditions, the relationship is easily understood:
the moment maps are related by a
\it projection\rm\ from the Lie algebra of the full group onto that of the
subgroup, while the notions of stability --- which are both formulated
as a set of algebraic conditions on subobjects --- are related by a
\it restriction\rm\ of the conditions to special subobjects determined by
the
subgroup.
In sections \ref{sect:StrSub} and \ref{sect:examples} we discuss some
specific types of subgroups of the gauge group and consider several
special cases of principal pairs. We show how our enlarged principal
pair framework encompasses them all and allows us to recover many
known Hitchin--Kobayashi correspondences.

A modification of the Hitchin--Kobayashi correspondences similar to
the one described in this paper, but more limited in scope, can be
found in \cite{OT}. The setting considered there corresponds only to
the special case discussed in Section \ref{sect:normal}.

\subparagraph{\bf Acknowledgements.}  We thank the mathematics departments
of the University of Illinois at Urbana-Champaign, the Universidad
Aut{\'o}noma de Madrid, the Mathematical Sciences Research Institute
of Berkeley and the Mathematical Institute of the University of
Oxford for their hospitality during various stages of this research.

\section{Principal pairs framework}
\subsection{Geometric Setting}\label{subs:geomsetting}
We summarize from [M] the structures we need for the definition and
analysis of principal pairs.
\subsubsection{} Let $K$ be a compact connected Lie group, and let
$G=K^{\CC}$ be its
complexification. Let $\klie$ (resp. $\glie$) be the Lie algebra of
$K$ (resp. $G$).
\subsubsection{} Fix a faithful unitary representation $\rho_a:K\to
\U(\WW_a)$, where
$\WW_a$ is a finite dimensional Hermitian vector space. Denote as
well by $\rho_a$ the induced holomorphic representation of $G$ on
$\GL(\WW_a)$ and let also $\rho_a:\glie\to\End(\WW_a)$ be the induced
Lie algebra representation. Let $*:G\to G$ be the Cartan involution.
Use this to define an Hermitian form on $\glie$, by:
$$\la u,v\ra = \Tr(\rho_a(u)\rho_a(v)^*).$$
Let $|u|=\la u,u\ra^{1/2}$. The restriction of $\la\ ,\ \ra$ to
$\klie$ will be used to identify ($K$-equivariantly) $\klie\simeq\klie^*$.
\subsubsection{} Let $X$ be a compact Kaehler manifold of complex dimension
$n$. Let
$\omega$ be the symplectic form of $X$, and denote as usual
$\omega^{[k]}=\omega^k/k!$. Unless otherwise stated, in the integrals
of functions on $X$ the volume form $\omega^{[n]}$ will be used. Let
$\Lambda:\Omega^{*+2}(X)\to\Omega^*(X)$ denote the adjoint of wedging
with $\omega$.
\subsubsection{} Let $P_K\to X$ be a $K$-principal bundle, and let
$P_G=P_K\times_K G$
be the $G$-principal bundle associated to $P_K$. Let
$\GGG=\Gamma(P_K\times_{\Ad}K)$ be the gauge group of $P_K$, and similarly
let $\GGC$ be the gauge group of $P_G$. Let
$\ad P_K=P_K\times_{\Ad}\klie$ and
${\ad}P_G=P_G\times_{\Ad}\glie$. Then $\Lie\GGG=\Omega^0({\ad}P_K)$ and
$\Lie\GGC=\Omega^0({\ad}P_G)$. Use $\la\ ,\ \ra$ to define an $L^2$ inner
product on $\Omega^0({\ad}P_K)$, and use this to obtain an inclusion
$\Lie\GGG\subset(\Lie\GGG)^*$.
\subsubsection{} Let $\Y$ be a complete Kaehler manifold with a Hamiltonian
action of
$K$, which we denote by $\rho:K\times \Y\to \Y$, and let
$$\mu:\Y\to\klie^*$$
\noindent be the corresponding moment map. Suppose that the action of $K$ on
$\Y$ respects the complex structure,
$I_{\Y}\in\End(T\Y)$, of $\Y$. Then, since $I_{\Y}$ is integrable,
there is a unique extension of the $K$-action on $\Y$ to a unique
holomorphic action of $G$.
\subsubsection{} Let $\YY_K=P_K\times_K \Y$ be the bundle with fibre $\Y$
associated
to $P_K$. Using the extension of the $K$-action to $G$, we can think
of $\YY_K$ as a $\Y$-bundle associated to $P_G$. In that case we
write $\YY_G=P_G\times_G \Y$. If the distinction is not crucial, we
write simply $\YY$.
\begin{definition}
We will say that the data $(P_K, \Y,\rho)$ determines a {\bf
symplectic pair type}. We will abbreviate this to $(P_K,\YY_K)$ where
$\YY_K=P_K\times_{\rho}\Y$. There is a corresponding {\bf complex pair type}
determined by the data $(P_G, \Y,\rho)$\ (abbreviated to
$(P_G,\YY_G)$ where $\YY_G=P_G\times_{\rho}\Y$).
\end{definition}
\subsection{The configuration space of pairs}
Let ${\AAA}$ be the space of connections on $P_K$ and let $\CCC$ be
the space of $G$-invariant almost complex structures on $P_G$ for
which the projection $P_G\to X$ is a pseudoholomorphic map.
\begin{prop} (see \S 2.2 in \cite{M})
There is a bijective correspondence between $\AAA$ and $\CCC$.
\end{prop}
The group
$\GGG$ acts on $\AAA$, and the group $\GGC$ acts on $\CCC$. In both cases
the
action is by pullback.  Thus we get a natural action of $\GGC$ on
$\AAA$. This action extends that of $\GGG$. Finally, let
$\AAA^{1,1}\subset\AAA$ be the set of connections whose curvature is of 
type $(1,1)$. The set $\AAA^{1,1}$ 
corresponds to
integrable complex structures on $P_G$, and is clearly
$\GGC$-invariant.  The group $\GGG$ acts on the space of
sections $\Gamma(\YY)$
and, since the action of $K$ on $\Y$ extends to an action of $G$, we
deduce that $\GGC$ acts on $\Gamma(\YY)$ extending the action of
$\GGG$.
We now define covariant derivations of sections of $Y=Y_K=Y_G$,
generalizing the usual definition in the vector bundle case (i.e.,
when $\Y$ is a vector space and the action of $K$ on $\Y$ is linear).
Any connection
$A\in\AAA$ induces a projection
$$\alpha_A:T\YY\to T\YY_v=\Ker d\pi,$$
where $\pi:\YY\to X$ is the projection ($T\YY_v$ is then the vertical
tangent bundle). Given any $\Phi\in\Omega^0$ we define
$$d_A\Phi:=\alpha_A\circ d\Phi\in\Omega^1(X;\Phi^*T\YY_v).$$
To generalize the notion of $\ov{\partial}$ operators on vector
bundles, observe that the complex structure on $\Y$ induces a
structure of complex vector bundle on $T\YY_v\to \YY$. Hence we can
decompose
$$d_A\Phi=\partial_A\Phi+\ov{\partial}_A\Phi,$$
where $\partial_A\in\Omega^{1,0}(X;\Phi^*T\YY_v)$ and
$\ov{\partial}_A\Phi\in\Omega^{0,1}(X;\Phi^*T\YY_v)$.
\begin{remark} If $A$ is an integrable connection on $P_K$, then a section
$\Phi$ satisfies the condition $\ov{\partial}_A\Phi=0$ if and only if
it is holomorphic with respect to the holomorphic structure induced
by $A$ on $\YY_G$.
\end{remark}
\begin{definition}\hfil
\begin{itemize}
\item A {\bf symplectic pair on $(P_K,\YY_K)$} consists of the pair
$(A,\Phi)$, where $A$\ is an integrable connection on $P_K$ and $\Phi$ is a
section of $\YY_K$ such that $\ov{\partial}_A\Phi=0$.
\item A {\bf complex pair on $(P_G,\YY_G)$} consists of the pair
$(I_G,\Phi)$, where $I_G\in\CCC$\ is a $G$-invariant
holomorphic structure on $P_G$ and $\Phi$ is a holomorphic section of
$\YY_G$.
\end{itemize}
\end{definition}
\begin{definition}  The {\bf configuration space of symplectic
pairs on $(P_K,\YY)$}, denoted by $\XXX(P_K,\YY_K)$ or simply
$\XXX_K$, is the subspace of $ \AAA^{1,1}\times\Gamma(\YY_K)$
defined by the condition that $\ov{\partial}_A\Phi=0$. Similarly, the
{\bf configuration space of complex pairs on $(P_G,\YY_G)$}, denoted
by $\XXX(P_G,\YY_G)$ or simply $\XXX_G$, is the subspace of
$\CCC\times\Gamma(\YY_G)$ defined by the condition that $\Phi$ is
holomorphic with respect to
$I_G$.
\end{definition}

Since $\YY_K=\YY_G$, there is a $\GGC$-equivariant bijection
between these configuration spaces.

\subsection{Moment map equations}\label{ss:mmeq}
It is an observation which goes back to the work of Atiyah and Bott
\cite{AB} and Donaldson \cite{D}that the set of connections
$\AAA$ carries a natural
symplectic structure, which can be defined by combining the Kaehler
structure on $X$ with the bi-invariant metric on $\klie$. Taking the
restriction of this symplectic structure we get a symplectic
structure on the smooth locus of the set $\AAA^{1,1}$. On the other
hand it is possible to define a symplectic structure on
$\Gamma(\YY)$ using the symplectic structure on
$\Y$ (see \S 4.2 in \cite{M} for details).
Combining both structures we get a symplectic structure on the smooth
locus of $\AAA^{1,1}\times\Gamma(\YY)$, whose restriction to
$\XXX(P_K,\YY_K)$ is also symplectic. It turns out that the action of
$\GGG$ on $\XXX(P_K,\YY_K)$ preserves this structure and is
Hamiltonian. Its moment map
$$\mu_{\GGG}:\AAA^{1,1}\times\Gamma(\YY)\to
\Lie \GGG\subset (\Lie\GGG)^*$$
can be computed to be
\begin{equation}
\mu_{\GGG}(A,\Phi)=\Lambda F_A+\mu(\Phi).
\label{eqn:defxi}
\end{equation}
\noi This formula should be interpreted as follows:
$\Lambda F_A$ is a section of ${\ad}P_K$, and
using the $K$-equivariant isomorphism $\klie\simeq\klie^*$ we can
also regard $\mu(\Phi)$ as a section of ${\ad}P_K$. Finally, the
$L^2$ inner product in $\Omega^0({\ad}P_K)$ gives the inclusion
$\Omega^0({\ad}P_K)\subset\Omega^0({\ad}P_K)^*=(\Lie\GGG)^*$.
In general the moment map of a hamiltonian action is unique only
up to addition of central elements in the (dual of the) Lie algebra.
In particular, if $c\in\klie$ is any
central element, then we can take the moment map to be
\begin{equation}
\mu_{\GGG,c}(A,\Phi)=\Lambda F_A+\mu(\Phi)-c.
\label{eqn:defxic}
\end{equation}
\subsection{Maximal weights and definition of stability}
\label{pesmaximal}
\begin{definition}\label{defn:Simplepairs}
We will say that an element $s\in\Omega^0({\ad}P_G)$ is semisimple if
for any $x\in X$, taking a
$G$-equivariant identification $({\ad}P_G)_x\simeq\glie$,
$s(x)\in\glie$ is semisimple.
We will say that a pair $(A,\Phi)\in\XXX(P_K,E)$ is {\bf simple} if
there is no semisimple element in $\Omega^0({\ad}P_G)$ which leaves
$(A,\Phi)$ fixed by the infinitesimal action. Observe that if a pair
$(A,\Phi)$ is simple, then any pair in the orbit
$\GGC(A,\Phi)$ is also simple.
\end{definition}
Let $\W=P_K\times_{\rho_a}\WW_a$. Any connection $A\in\AAA^{1,1}$
induces a
$\ov{\partial}$-operator
$$\ov{\partial}_A:\Omega^0(\W)\to\Omega^{0,1}(\W)\ .$$
\begin{definition}
\label{defn:degchi}
Let $\chi\in\Omega^0({\ad}P_K)$. We will say that $\chi$ {\bf induces an
$A$-holomorphic filtration} if the following two conditions
are satisfied:
(1) all the eigenvalues of $\rho_a(\imag\chi)$ acting on $\W$ are
constant;
(2) if $\alpha_1<\dots<\alpha_r\in\RR$ are the different eigenvalues
of
$\rho_a(\imag\chi)$, and we define
$$\W^k=\bigoplus_{j\leq k}\Ker(\alpha_j\Id-\rho_a(\chi))\subset \W
\quad ;\quad 1\leq k\leq r$$
\noi then the filtration
$\W^1\subset \W^2\subset\dots\subset \W^r=\W$
is holomorphic w.r.t. $A$, that is,
$$\ov{\partial}_A(\W^k)\subset \Omega^{0,1}(\W^k)\ .$$
\noi When these conditions hold, we define the {\bf degree of $\chi$} to be
$$\deg(\chi)=\alpha_r\deg(\W)+
\sum_{k=1}^{r-1}(\alpha_k-\alpha_{k+1})\deg(\W^k),$$
where for any vector bundle $\W'$ we denote
$\deg(\W')={2\pi}\la c_1(\W')\cup [\omega^{[n-1]}],[X]\ra.$
(Here $[\omega^{[n-1]}]$ denotes the cohomology class represented by
the form $\omega^{[n-1]}$ and $[X]\in H_{2n}(X;\ZZ)$ is the
fundamental class of $X$.)
\end{definition}
The condition that $\chi\in\Omega^0({\ad}P_K)$ induces
$A$-holomorphic filtration is in fact independent of the chosen
representation
$\rho_a$ (provided it is faithful). Indeed, the induced
holomorphic filtrations are in bijection with holomorphic reductions
of the structure group of $P_G$ to parabolic subgroups of
$G$ (see \S 2.7 in \cite{M}).
Finally, observe that the preceding definitions make sense when the
section $\chi$ is defined only on $X-X_0$ where $X_0$ is a closed
complex subset of $X$ and $X-X_0$ has codimension at least two (in
this situation $\chi$ defines a filtration of $\W|_{X_0}$, and the
degree of each subbundle in the filtration is well defined, see
\cite{UY}). This will be relevant below, when we will consider
elements of
$\hlie(X_0)$.
\begin{definition}\label{defn:maxwghts} Let $(\Y,\omega_{\Y},I_{\Y})$ be a
Kaehler manifold with a Hamiltonian action of a compact Lie group
$K$, and corresponding moment map $\mu: \Y\rightarrow\klie^*$. Let
$x\in \Y$ be any point, and take an element
$s\in\klie$.
We define the {\bf maximal weight $\lambda(x;s)$ of the action of $s$
on $x$} to be
$$\lambda(x;s)=\lim_{t\to \infty}\la\mu(e^{\imag ts}x),s\ra
\in\RR\cup\{\infty\}.$$
\end{definition}
Notice that the maximal weight is $K$-invariant in the sense that
$\lambda(x;s)=\lambda(kx;ksk^{-1})$ for any $k\in K$ (this follows
from equivariance of the moment map). For any $\Phi\in\Gamma(\YY)$
and $\chi\in\Omega^0({\ad}P_K)$ we can thus define a map
$\lambda(\Phi;\chi):X\to\RR\cup\{\infty\}$ by using any local frame
to identify $Y_x=\Y$ and $({\ad}P_K)_x=\klie$, and setting
$$\lambda(\Phi;\chi)(x)=\lambda(\Phi(x);\chi(x))\ .$$
\begin{definition}\label{defn:lambdac}Let $(A,\Phi)$ be any pair in
$\AAA^{1,1}\times\Gamma(\YY)$. Let $\chi$ be a section in
$\Omega^0({\ad}P_K)$ which induces a $A$-holomorphic
filtration and let $c$ be a central element in $\klie$. Define the
{\bf maximal weight}
\begin{equation}\label{eqtn:maxweight}
\lambda_c(A,\Phi;\chi)=
\deg(\chi)+\int_{X}\lambda(\Phi;\chi)-\int_{X}\la\chi,c\ra.
\end{equation}
\end{definition}
\begin{remark} We can identify $\deg(\chi)$
as a maximal weight for the action of $\chi$ on
$\AAA^{1,1}$ (cf. Lemmas 4.2 and 4.3 in \cite{M}).
Thus $\lambda_c(A,\Phi;\chi)$
is in fact the maximal weight in the sense of Definition
\ref{defn:maxwghts} for the action of $\chi$ on
$\AAA^{1,1}\times\Gamma(\YY)$,
and using the moment map $\mu_{\GGG,c}$ for product
$\GGG$-action on $\AAA^{1,1}\times\Omega^0(\YY)$.
\end{remark}
\begin{definition}\label{defn:Stability}
Let $(A,\Phi)\in\XXX(P_K,\YY)$, and let
$c$ be a central element in $\klie$.
We say that $(A,\Phi)$ is {\bf $c$-stable} if for
\begin{itemize}
\item any open subset $X_0\subset X$ whose complementary
has complex codimension $\geq 2$ and
\item any $\chi\in\Omega^0(X_0;{\ad}P_K)$ which induces a $A$-holomorphic
filtration, we have
\begin{equation}\label{eqtn:chi0}
\deg(\chi)+\int_{X_0}\lambda(\Phi;\chi)-\int_{X_0}\la\chi,c\ra>0
\end{equation}
(if the integral $\int_{X_0}\lambda(\Phi;\chi)$ is equal to
$\infty$ then, since the other terms are finite numbers, we consider the
left hand side as being greater than zero).
\end{itemize}
\end{definition}

\begin{remark} Notice that if $X_0=X$, then (\ref{eqtn:chi0}) says
$\lambda_c(A,\Phi;\chi)>0$.
\label{rk:maxweight}
\end{remark}

\subsection{A Hitchin--Kobayashi correspondence}
The main result in \cite{M} describes which orbits of the $\GGC$
action on $\XXX(P_K,\YY)$ contain solutions to the moment map
equation $\mu_{\GGG}(A,\Phi)=c$.

\begin{theorem}\label{theorem:HK}
(Hitchin--Kobayashi Correspondence for principal pairs) Let
$c\in\klie$ be a central element. Let $(A,\Phi)\in\XXX(P_K,\YY)$ be a
simple pair. Then $(A,\Phi)$ is $c$-stable if and only if there
exists a gauge transformation
$g\in\GGC$ such that $g(A,\Phi)=(g^*A,g(\Phi))$
satisfies the generalized vortex equation
\begin{equation}
\Lambda F_{g^*A}+\mu(g(\Phi))=c
\label{eqn:vortex}
\end{equation}
\noi Furthermore, if $(A,\Phi)$ is $c$-stable then any two such gauge
transformations $g,g'\in\GGC$ satisfy $g'g^{-1}\in\GGG$.
\end{theorem}
\begin{remark} The equation (\ref{eqn:vortex}) generalizes the vortex
equations, which arise in the case $\V=\CC^n$ and $K=\U(n)$ acting
through the fundamental representation. The vector bundle case (see
Section \ref{sss:vectorbdl}) of this result was proved by Banfield in
\cite{Ba}.
\end{remark}
\subsection{The vector bundle case}\label{sss:vectorbdl}
Suppose that $\Y=\V$ is a Hermitian vector space with the Kaehler
structure given by the Hermitian metric, and that the action
$\rho:K\times \V\to \V$ is linear, so that it comes from a morphism of
groups $\rho_0:K\to \U(\V)$. Then $\VV=P\times_{\rho}\V$ is a vector
bundle, and the stability condition can be greatly simplified.
To explain this we need to make a definition. Let $A\in\AAA^{1,1}$,
and suppose that $\chi\in\Omega^0(X_0;{\ad}P_K)$ induces an
$A$-holomorphic filtration on $W=P_K\times{\rho_a}\W_a$,
(as in Definition \ref{defn:degchi}) with eigenvalues
$\imag\{\alpha_1\le\alpha_2\le\dots\le\alpha_r\}$ (here
$\alpha_j\in\RR$).  Then $\imag\rho_0(\chi)$
induces an $A$-holomorphic filtration on $V$, also consisting of
subbundles spanned by eigenvectors.
The eigenvalues will in general be combinations
(which depend on the representation $\rho_0$) of the $\alpha_i$
(see section \ref{sect:examples} for examples).
%
%

\begin{definition} Let $\VV^-(\chi)\subset\VV$ denote
the subbundle spanned by the eigenvectors of the endomorphism
$\imag\rho_0(\chi)$ with non-positive eigenvalue.
By assumption $\VV^-(\chi)$ is a holomorphic subbundle.
\end{definition}
\begin{lemma}
A pair $(A,\Phi)\in\XXX(P_K,\YY)$ is $c$-stable if and only if: for
any
$\chi\in\Omega^0(X_0;{\ad}P_K)$ inducing an $A$-holomorphic filtration
and such that $\Phi\subset \VV^-(\chi)$ we have
$$\deg(\chi)-\int_{X_0}\la\chi,c\ra > 0.$$
\end{lemma}
\begin{pf} Let $s\in\klie$, and let $\V^-(s)\subset \V$ be the subspace
spanned by the eigenvectors of the endomorphism $\imag \rho_0(s)$
with eigenvalue $\leq 0$. Then for any $x\in \V$ we have
$$\lambda(x;s)=\left\{\begin{array}{ll}
0 &  \text{if $x\in \V^-(s)$} \\
\infty & \text{if $x\notin \V^-(s)$.}
\end{array}\right.$$
This proves the result.
\end{pf}

This result relates our general stability condition with the notion
of stability given by Banfield in \cite{Ba} in the case of vector
bundle pairs.

\section{Subgroups of the gauge group}
\subsection{The subgroup setting}\label{subs:subgp}
Let $\HHH\subset\GGG$ and $\HHC\subset\GGC$ be Lie subgroups (with
respect to the $C^{\infty}$ topology), and consider their respective
Lie algebras
\begin{align*}
\hlie & = \{ s\in\Omega^0({\ad}P_K) \mid \exp(ts)\in\HHH\ \forall t\in\RR\}, \\
\hlie^{\CC} & = \{ s\in\Omega^0({\ad}P_G) \mid
\exp(ts)\in\HHC\ \forall t\in\RR\}.
\end{align*}
\begin{assumption}\label{ass:sub}
We make the following assumptions:
\begin{enumerate}
\item $\hlie^{\CC}=\hlie\otimes_{\RR}\CC$,
\item the map $\exp:\imag\hlie\to\HHC$ induces an isomorphism
$\imag\hlie\simeq\HHC/\HHH$.
\item  $\Omega^0({\ad}P_K)=\hlie^{\perp}\oplus\hlie$ is a splitting
of Fr\'echet spaces, with $\hlie^{\perp}$ orthogonal to $\hlie$ with
respect to the pairing $\int_X\la\ ,\ \ra:\Omega^0({\ad}P_K)\otimes
\Omega^0({\ad}P_K) \to\RR$.
\end{enumerate}
\end{assumption}
\begin{remark}  The first two conditions can be rephrased by saying that
$\HHC$ is the complexification of $\HHH$.
\end{remark}
\begin{definition}\label{defn:subgp}Given a subgroup $\HHC$, we can define
$\HHC$-invariant
subsets $\XXX_{\HHH}\subset\XXX(P_K,\YY)$. A pair
$(A,\Phi)\in\XXX_{\HHH}$ will be called {\bf an $\HHH$-pair of type
$(P_K,\YY)$}.
We say that $(\HHH,\HHC,\XXX_{\HHH})$ defines a {\bf subgroup
setting} if
\begin{itemize}
\item the subgroups $\HHH\subset\GGG$ and $\HHC\subset\GGC$ satisfy
Assumptions \ref{ass:sub},
\item $\XXX_{\HHH}\subset\XXX(P_K,\YY)$ is an $\HHC$-invariant subset, and
\item $d_A(\hlie)\subset \hlie^1$ where $A$ is any connection in an
$\HHH$-pair of type $(P_K,\YY)$ and
$$\hlie^1=\Omega^1(X)\wo \hlie\subset\Omega^1({\ad}P_K)\ .$$
\end{itemize}
\end{definition}
\noi We will see later in Lemma \ref{condinv} that if a connection $A$
satisfies the condition $d_A(\hlie)\subset \hlie^1$ then any other
connection in the orbit $\GGC A$ also satisfies it. For future
reference, we make two more definitions:
\begin{definition}
Denote by
$\pi_{\hlie}:\Omega^0({\ad}P_K)\to\hlie$ the projection induced by the
splitting $\Omega^0({\ad}P_K)=\hlie^{\perp}\oplus\hlie$ in (3).
\end{definition}

\begin{definition}
If $X_0\subset X$ is an open subset such that $X-X_0$ is a complex
subset of codimension at least two, define
$$\hlie(X_0)=\left\{ \sigma\in\Omega^0(X_0;{\ad}P_K)
\mid \forall s\in\hlie^{\perp},\ \int_{X_0}\la\sigma,s\ra=0\right\}.$$
\label{def:restriction}
\end{definition}

\begin{remark} The need to formulate
Definition \ref{def:restriction} in this way can be understood as
follows. Recall that the definition of $c$-stability (see Definition
\ref{defn:Stability}) uses the positivity of a certain integral
defined over open subsets $X_0\subset X$ where
$\codim X\setminus X_0\geq 2$. The integral can be interpreted as the
maximal
weight of the pair $(A,\Phi)$ with respect to an element in the Lie
algebra for the gauge group of the restriction of $P_K$ to
$X_0$ (see Remark \ref{rk:maxweight}). Hence if $\dim X\geq 2$,
it is not enough to check positivity of the maximal weights with
respect only to actual elements of the Lie algebra of the gauge
group\footnote{In the Hitchin--Kobayashi correspondence for vector
bundles without additional structure this translates into the need to
consider reflexive subsheaves.}.

We might describe the situation as follows:
consider on $X$ the topology whose open sets are of the form
$X_0\subset X$ with $\codim X\setminus X_0\geq 2$, and let
$\GGGS$ be the sheaf such that for any open set $X_0\subset X$
the sections $\Gamma(\GGGS;X_0)$ are the automorphisms of
$P_K|_{X_0}$. Then $\GGGS$ is a sheaf of groups and, strictly
speaking, the Hitchin--Kobayashi correspondence should be understood
as a version of the Kempf--Ness theorem for $\GGGS$ (acting on the
sheaf defined by considering connections and sections defined over
open subsets). It is only when $\dim X=1$ that one can identify the
sheaf $\GGGS$ with the gauge group $\GGG$ and properly say that the
Hitchin--Kobayashi correspondence is a theorem {\it \`a la}
Kempf--Ness for the infinite dimensional group $\GGG$.

{}From this point of view it is clear that rather than a subgroup of
the gauge group what we need is to specify a subsheaf of $\GGGS$.
Given a subgroup $\HHH\subset\GGG$, Definition \ref{def:restriction}
is intended to provide a subsheaf $\HHHS\subset\GGGS$ (although this
is done, strictly speaking, at the level of Lie algebras), in such a
way that the main result of this paper (Theorem \ref{theorem:main})
is a theorem {\it \`a la} Kempf--Ness for the sheaf $\HHHS$.
\label{rk:restriction}
\end{remark}

\subsection{The vortex equations for $\HHH$}
\label{ss:mmeqH}
The restriction of the Hamiltonian action of $\GGG$ on
$\XXX(P_K,\YY_K)$ to $\HHH\subset\GGG$ is also Hamiltonian. Thus we can
define a moment map for the $\HHH$-action on the smooth locus of
$\XXX(P_K,\YY_K)$. This moment map for
$\HHH$ is the composition of
$\mu_{\GGG}:\XXX(P_K,\YY)\to\Lie \GGG \subset(\Lie\GGG)^*$ with
the projection $(\Lie \GGG)^*\to \hlie^*$ induced by the
inclusion. Equivalently, it is the composition of the map
$\mu_{\GGG}:\XXX(P_K,\YY)\to\Lie \GGG$ with the orthogonal
projection $\pi_{\hlie}:\Lie\GGG\to\hlie$ and the inclusion
$\hlie\subset\hlie^*$ induced by restricting the $L^2$ inner
product on $\Omega^0({\ad}P_K)$ to $\hlie$.
On the other hand, the restriction of the symplectic structure on
$\XXX(P_K,\YY)$ gives a symplectic structure on $\XXX_{\HHH}$.
Since $\XXX_{\HHH}$ is $\HHH$-invariant, we deduce that $\XXX_{\HHH}$
carries a Hamiltonian action of $\HHH$. Its moment map, as described
above, is thus
$$\mu_{\HHH}(A,\Phi)=\pi_{\hlie}(\Lambda F_A+\mu(\Phi))\ .$$

\begin{definition}
Let $c$ be a constant central element in $\klie$ and let
$c_{\HHH}$ be its projection onto $\hlie$.
We say an $\HHH$-pair of type $(P_K,\YY)$ satisfies the {\bf
$(\HHH,c_{\HHH})$-vortex equations} if there exists $h\in\HHC$ such that
\begin{equation}\label{eqtn:HcVortex}
\pi_{\hlie}(\Lambda F_{h(A)}+\mu(h(\Phi)))=c_{\HHH}
\end{equation}
\end{definition}
\subsection{Stability for elements of $\XXX_{\HHH}$}
\begin{definition}\label{defn:HSimplepairs}
We will say that an element $s\in\hlie^{\CC}$ is semisimple if, seen
as a section in $\Omega^0({\ad}P_G)$, it is semisimple (see Definition
\ref{defn:Simplepairs}).
We will say that a pair $(A,\Phi)\in\XXX_{\HHH}$ is {\bf simple} if
there is no semisimple element in $\hlie^{\CC}$ which leaves
$(A,\Phi)$ fixed by the infinitesimal action.
Observe that if a pair $(A,\Phi)$ is simple, then any pair in the
orbit $\HHC(A,\Phi)$ is also simple.
\end{definition}
\begin{definition}
\label{defn:HcStability}
Let $(A,\Phi)\in\XXX_\HHH$, and let
$c_{\HHH}$ be a constant central element in $\hlie$.
We will say that $(A,\Phi)$ is {\bf
$(\HHH,c_{\HHH})$-stable} if
$$\deg(\chi)+\int_{X_0}\lambda(\Phi;\chi)-\int_{X_0}\la\chi,c_{\HHH}\ra>0$$
\noi whenever
\begin{itemize}
\item  $X_0\subset X$ is an open set such that $X-X_0$ is a complex
subset of codimension at least two, and
\item $\chi\in\hlie(X_0)$ induces an $A$-holomorphic filtration.
\end{itemize}
(if the integral $\int_{X_0}\lambda(\Phi;\chi)$ is equal to
$\infty$ then, since the other terms are finite numbers, we consider the
left hand side as being to be greater than zero).
\end{definition}
\begin{remark} (The vector bundle case) When we are in the situation
described in \ref{sss:vectorbdl}, the definition of
$(\HHH,c_{\HHH})$-stability simplifies exactly in the same
way as in the definition of $c$-stability.
\end{remark}
\section{The Main Theorem}\label{sect:main}
Our Main Theorem provides the Hitchin--Kobayashi Correspondence for
$\HHH$-pairs of type $(P_K,\YY)$.
\begin{theorem}
\label{theorem:main}{\bf (Main Theorem)}
Let $(\HHH,\HHC,\XXX_{\HHH})$ define a subgroup setting, as in
Definition \ref{defn:subgp}. Let
$c_{\HHH}$ be a constant central element in $\hlie$. Let
$(A,\Phi)\in\XXX_{\HHH}$ be a simple pair.
Then $(A,\Phi)$ is $(\HHH,c_{\HHH})$-stable if and only if there
exists
$h\in\HHC$ such that
\begin{equation}
\pi_{\hlie}(\Lambda F_{h(A)}+\mu(h(\Phi)))=c_{\HHH}.
\label{lequacio}
\end{equation}
Furthermore, if two different $h,h'\in\HHC$ solve equation
(\ref{lequacio}), then there exists $k\in\HHH$ such that
$h'=kh$.
\end{theorem}
The proof will follow very closely that of Theorem 2.19 in
\cite{M} (which in turn relies on \cite{Br}).
Consequently we will sketch the main steps and will only explain in
detail the new features of the proof. Note that there are many papers
proving similar results, and we refer to \cite{M} for a partial list
of them.
\subsection{Preliminaries}
To begin with, it is convenient to complete our spaces by means of
suitable Sobolev norms. Take $p>2n$. We use the $L^p_2$ norm on
$\GGC$ and the $L^p_1$ norm on $\Gamma(\YY)$. Let $A_0\in\AAA$ be any
smooth connection and write
$$\AAA_{L^p_1}=A_0+L^p_1(T^*X\otimes {\ad}P_K)\ .$$
\noindent With this choices, $(\GGC)_{L^p_2}$ is a Hilbert Lie
group which acts smoothly on $\AAA_{L^p_1}$ and on
$\Gamma(\YY)_{L^p_1}$. We have
$\Lie (\GGG)_{L^p_2}=L^p_2({\ad}P_K)$, and
$(\GGC)_{L^p_2}$ is the complexification of $(\GGG)_{L^p_2}$.
We also take the completions $\HHH_{L^p_2}$ and $\HHH^{\CC}_{L^p_2}$
of $\HHH$ and $\HHC$ with respect to the $L^p_2$ norm, which are Lie
subgroups of
$(\GGG)_{L^p_2}$ and $(\GGC)_{L^p_2}$ respectively.
The Lie algebra of $\HHH_{L^p_2}$ (resp. $\HHH^{\CC}_{L^p_2}$) is the
completion $\hlie_{L^p_2}$ (resp. $\hlie^{\CC}_{L^p_2}$) of $\hlie$
(resp. $\hlie^{\CC}$) with respect to the $L^p_2$ norm. Observe that
$\HHH_{L^p_2}$ (resp.
$\HHH^{\CC}_{L^p_2}$, $\hlie_{L^p_2}$, $\hlie^{\CC}_{L^p_2}$)
is the closure of $\HHH$ (resp. $\HHC$, $\hlie$, $\hlie^{\CC}$) in
$(\GGG)_{L^p_2}$ (resp. $(\GGC)_{L^p_2}$, $L^p_2({\ad}P_K)$,
$L^p_2({\ad}P_G)$).
Since the splitting $\Omega^0({\ad}P_K)=\hlie^{\perp}\oplus\hlie$ is
assumed to be of Fr\'echet spaces, it follows that $\pi_{\hlie}$
extends to a continuous operator
$$(\pi_{\hlie})_{L^p_2}:L^p_2({\ad}P_K)\to\hlie_{L^p_2}.$$
Most of the time we will  avoid writing the
Sobolev subscripts, and the Sobolev norms will be implicitly assumed.
Recall (see \ref{ss:mmeqH}) that $\XXX_\HHH$ admits on its smooth
locus a Kaehler structure which is preserved by the action of $\HHH$.
Furthermore, the action of
$\HHH$ admits a moment map $\mu_\HHH:\XXX^0_\HHH\to\hlie^*\cong\hlie$,
where $\XXX^0_{\HHH}$ is the smooth locus in $\XXX_{\HHH}$. It
follows from the general results in
\cite{M} that for any central $c_{\HHH}\in\hlie$ there exists a function
$\Psi:\XXX_\HHH\times\HHC\to\RR$, called the integral of the moment
map, which satisfies the following key properties:
\begin{enumerate}
\item If $s\in\hlie$, $\Psi((A,\Phi),e^{\imag s})
=\int_0^1\la\mu_{\HHH}(e^{\imag ts}(A,\Phi))-c_{\HHH},s\ra dt$, and
if $h\in\HHH$, then $\Psi((A,\Phi),ke^{is})=\Psi((A,\Phi),e^{is})$.
\item If $g,h\in\HHC$, then
$\Psi((A,\Phi),g)+\Psi(g(A,\Phi),h)=\Psi((A,\Phi),hg)$.
\item Let $\Psi_{(A,\Phi)}:\HHC\to\RR$ be the restriction of $\Psi$ to
$(A,\Phi)\times\HHC$. The element $g\in\HHC$ is a critical point
of $\Psi_{(A,\Phi)}$ if and only if $(A',\Phi')=g(A,\Phi)$ satisfies
$\Lambda F_{A'}+\mu(\Phi')=c_{\HHH}$.
\end{enumerate}
\noindent Note that (1) defines $\Psi$ and that (3) follows from (1) and
(2).
In view of the Kaehler structure on $\XXX_\HHH$ and the action of
$\HHH$, it makes sense to define maximal weights in the same way as we
did for the action of $\GGG$ on $\XXX$ in \ref{pesmaximal}. Thus,
as in Definition \ref{defn:lambdac}, we define
$$\lambda_{c_{\HHH}}(A,\Phi;\chi)=
\deg(\chi)+\int_{X}\lambda(\Phi;\chi)-\int_{X}\la\chi,c_{\HHH}\ra\ .$$
Furthermore, this also makes sense in the singular locus of
$\XXX_\HHH$ because
$\XXX_\HHH\subset\AAA\times\Gamma(\YY)$ and the latter space is smooth
and Kaehler. By Lemma 4.3 in  \cite{M}, we have the following.
\begin{lemma}
Let $(A,\Phi)$ be a pair in $\XXX_\HHH$. Take
$\chi\in\hlie$ and $c_{\HHH}$ a constant central element in
$\Lie(\HHH)$. If the maximal weight
$\lambda_{c_{\HHH}}((A,\Phi);\chi)<\infty$, then $\chi$ induces an
$A$-holomorphic filtration and
$$\lambda_{c_{\HHH}}((A,\Phi);\chi)=
\deg(\chi)+\int_X\lambda(\Phi;\chi)-\int_{X}\la\chi,c_{\HHH}\ra\ .$$
\label{desmaxpes}
\end{lemma}
%
In view of the Kaehler interpretation of the equation
$\pi_{\hlie}(\Lambda F_A+\mu(\Phi))=c_{\HHH}$, the uniqueness claim of the
theorem
now follows from general results on convexity of the integral of the
moment map (see (3) of Proposition 3.3 and Theorem 5.4 in \cite{M}).
The next two lemmas show that the property $d_A\hlie\subset
\hlie^1$ depends only on the $\HHC$ orbit of $A$.
\begin{lemma}
Let $A\in\AAA$ be such that $d_A\hlie\subset \hlie^1$. Then for any
$g\in\HHC$ we have $d_{g(A)}\hlie\subset\hlie^1$.
\label{condinv}
\end{lemma}
\begin{pf}
Use the decomposition
$\Omega^1({\ad}P_G)=\Omega^{1,0}({\ad}P_G)\oplus\Omega^{0,1}({\ad}P_G)$ to
split the covariant derivative as $d_A=\partial_A+\ov{\partial}_A$.
Then  $d_{g(A)}=(g^*)^{-1}\partial_Ag^*+g\ov{\partial}_Ag^{-1}.$
{}From this it follows easily that
$d_{g(A)}(\hlie^{\CC})\subset\Omega^1(X)\wo \hlie^{\CC}$. But since
$d_{g(A)}$ is a $K$-connection (i.e., it is compatible with Cartan
involution) we deduce from the above inclusion that
$d_{g(A)}\hlie\subset\hlie^1$.
\end{pf}

\begin{lemma}
Suppose that $A\in\AAA$ satisfies $d_A\hlie\subset\hlie^1$. Then
$\Delta_A\hlie=d_A^*d_A\hlie\subset\hlie$.
\label{DeltaA}
\end{lemma}
\begin{pf}
It suffices to check that $d_A^*\hlie^1\subset\hlie$. Let
$(\hlie^{\perp})^1=\Omega^1(X)\wo \hlie^{\perp}$. Given
$a,b\in\Omega^0({\ad}P_K)$, integration by parts gives
$$\int_X\la d_Aa,b\ra=-\int_X\la a,d_Ab\ra,$$
from which we deduce that $d_A \hlie^{\perp}\subset
(\hlie^{\perp})^1$. Let now $h\in \hlie^1$ and $l\in \hlie^{\perp}$.
Then
$$\int_X\la d_A^*h,l\ra=\int_X\la h,d_Al\ra=0,$$
so that $d_A^*\hlie^1\subset (\hlie^{\perp})^{\perp}=\hlie$.
\end{pf}

\subsection{Stability implies existence of solutions}
\label{stabexist}
This section and the next one provide only a sketch of the proof,
with the emphasis on the modifications required in the proof for the
case $\HHH=\GGG$.  See
\cite{Br,M} for more details. Let $(A,\Phi)\in\XXX_\HHH$ be a simple
and $(\HHH,c_{\HHH})$-stable pair. Let us prove that there is a
solution of (\ref{lequacio}) in the
$\HHC$ orbit of $(A,\Phi)$.
Let $\Met=\imag\hlie_{L^p_2}$, and consider the functional
$$\begin{array}{rcl}
\Psi_{\hlie}:\Met &  \longrightarrow &  \RR \\
s &  \mapsto &  \Psi_{\hlie}((A,\Phi),e^s).
\end{array}$$
The elements $s\in\Met$ moving the pair $(A,\Phi)$ to a solution of
equation (\ref{lequacio}) are precisely the critical points of
$\Psi_{\hlie}$. For technical reasons it is convenient to restrict
ourselves to the set
$$\MetB=\{s\in\Met\mid
\|\pi_{\hlie}(\Lambda F_{e^sA}+\mu(e^s\Phi)-c)\|_{L^p}\leq B\},$$
for some $B\geq 0$. To do this safely it is necessary to check that
the critical points of the restriction of $\Psi_{\hlie}$ to $\MetB$
are also critical points of $\Met$. It is here that one needs the
pair $(A,\Phi)$ to be simple and the inclusion
$\Delta_{h(A)}\hlie\subset\hlie$ for any $h$ and
$A$ given by Lemmas \ref{condinv} and \ref{DeltaA} (see \S6.2.1 in
\cite{M}).
The key step of the proof is to deduce from
$(\HHH,c_{\HHH})$-stability the existence of positive constants
$C_1,C_2$ such that for any
$s\in\MetB$
\begin{equation}
|s|_{C^0}\leq C_1\Psi_{\hlie}(s)+C_2.
\label{desclau}
\end{equation}
This can be done following word by word \S6.2.1 in \cite{M} (and note
that it is here that one needs to use $A\in\AAA^{1,1}$, to invoke a
theorem of Uhlenbeck and Yau on weak $L^p_1$ bundles; see
\cite{Br,M,UY}). Finally, it follows from the inequality that there
exists some
$s\in\MetB$ which is a critical point of $\Psi_{\hlie}$, and hence gives a
solution to equation (\ref{lequacio}) (see \cite{Br, M}) (to prove
this one takes a minimizing sequence $\{s_j\}\subset\MetB$ and uses
the inequality to deduce convergence ---in a suitable sense--- of a
subsequence). Furthermore, one can prove that $s$ is smooth using
elliptic regularity (see \cite{Br}).
\subsection{Existence of solution implies stability}
Let $(A,\Phi)\in\XXX_\HHH$ be a simple pair, and assume that there
exists $h\in\HHC$ solving equation (\ref{lequacio}). Our aim is to
prove that $(A,\Phi)$ is $(\HHH,c_{\HHH})$-stable. In view of Lemma
\ref{desmaxpes}, this is equivalent to showing that for any open
$X_0\subset X$ whose complement in $X$ has complex codimension
$\geq 2$ and any $\chi\in\hlie(X_0)$ we have
$$\lambda_{c_{\HHH}}((A,\Phi);\chi)>0.$$
This is done in two steps. First one checks that
$\lambda_{c_{\HHH}}(h(A,\Phi);\chi)>0$ for any $\chi\in\hlie(X_0)$, i.e.,
that $h(A,\Phi)$ is indeed $(\HHH,c_{\HHH})$-stable. In particular,
this implies, using the same methods as in
\ref{stabexist}, that we have an inequality of the type
(\ref{desclau}). In the second step one uses this inequality to
relate the maximal weights at $h(A,\Phi)$ to those at
$(A,\Phi)$, thus deducing from the $(\HHH,c_{\HHH})$-stability of
$h(A,\Phi)$
that $(A,\Phi)$ has also to be $(\HHH,c_{\HHH})$-stable. See \S6.3 in
\cite{M} for more details.
\section{Stability Simplification Conditions (SSC)}\label{sect:SSC}
For some choices of the data
$(K,P_K,\Y,\rho)$ the general stability condition
\ref{defn:Stability} can be simplified. Classical examples
are the stability condition for holomorphic vector bundles and for
holomorphic pairs. In both cases the general stability condition,
which refers to any filtration of the vector bundle by holomorphic
subbundles (or in general reflexive subsheaves), can in fact be
reduced to the same condition considered only for subbundles. Our aim
in this section is to find a general condition (which we call
Stability Simplification Condition) which implies the possibility of
simplifying the stability condition. Our results are similar in
spirit to those by Schmitt in \cite{Sch}.
\subsection{Statement of the conditions}

\begin{definition}
\label{defn:SSC} Fix data $(K,P_K,\Y,\rho)$. For any open subset $U\subset
X$ and subset $S\subset\Omega^0({\ad}P_K)$ we denote by $S^{\perp}|_U$
the image of $S^{\perp}\subset\Omega^0(X;{\ad}P_K)$ under the
restriction map $\Omega^0(X;{\ad}P_K)\to \Omega^0(U;{\ad}P_K)$.

We say that the  Stability Simplification Condition (or SSC for
short) applies if there is a subset $S\subset\Omega^0({\ad}P_K)$ such
that for any $A\in\AAA^{1,1}$:
\begin{itemize}


\item (SSC1) for any $X_0\subset X$ whose complement has complex
codimension at least two, any $\sigma\in\Omega^0(X_0;{\ad}P_K)$ inducing
a $\ov{\partial}_A$ holomorphic filtration can be written
$$\sigma=s_1+\dots+s_r,$$
where each $s_i$ belongs to the orthogonal complement of
$S^{\perp}|_{X_0}$ and induces a $\ov{\partial}_A$ holomorphic filtration.

\item (SSC2) for any pair $(A,\Phi)\in\XXX(P_K,\YY)$ the condition
$\lambda(\Phi;\sigma)>0$ is equivalent to the condition
$\lambda(\Phi;s_i)>0$
for any $i$.

\end{itemize}
If we are in the situation described in \ref{sss:vectorbdl}, i.e. if
$\Y=\V$ is a vector space, then we can replace condition (SSC2)
by the following one.
\begin{itemize}
\item (SSC2') for any pair
$(A,\Phi)\in\XXX(P_K,\VV)$ the condition
$\Phi\subset \VV^-(\sigma)$ is equivalent to the condition
$\Phi\subset \VV^-(s_i)$ for any $i$.
\end{itemize}
\end{definition}

\begin{remark} The need to formulate (SSC1) in this way can be traced
back to the same considerations as those behind Definition
\ref{def:restriction} (see Remark \ref{rk:restriction})
\end{remark}
The following is almost a tautology.
\begin{lemma}
Let $(A,\Phi)\in\XXX(P_K,\YY)$ be any pair, and let $c\in\klie$ be a
central element. Then $(A,\Phi)$ is $c$-stable if and only if for any
$X_0\subset X$ whose complementary has complex codimension
$\geq 2$ and any $s\in S|_{X_0}$ inducing a $A$-holomorphic filtration we
have
$$\deg(s)+\int_{X_0}\lambda(\Phi;s)-\int_{X_0}\la s,c\ra >0.$$
\label{lemma:SSC}
\end{lemma}
\begin{remark} If $\Y=\V$ is a vector space, then there is an obvious
simplification of Lemma \ref{lemma:SSC} using the results in
\ref{sss:vectorbdl}. The details are left to the reader.
\end{remark}
The following lemma will be useful later to prove that SSC holds in
some situations.
\begin{lemma}
If SSC holds for $(K_i,P_{K_i},\rho_i,\Y)$ for $i=1,2$ and the
actions
$\rho_1,\rho_2$ of $K_1,K_2$ on $\Y$ commute, then
$\rho_1\times\rho_2$ defines an action of
$K_1\times K_2$ on $\Y$ and SSC holds for $(K_1\times K_2,P_1\times
P_2,\rho_1\times\rho_2,\Y)$.
\end{lemma}
\begin{pf} It follows from the fact that if $x\in \Y$ and $s_i\in K_i$
for $i=1,2$, then we have this equality between maximal weights:
$$\lambda(x;s_1+s_2)=\lambda(x;s_1)+\lambda(x;s_2).$$
\end{pf}

\begin{example}\label{ex:SSC1} Take $\V=\CC^n$, $K=\U(n)$ and $\rho$ the
fundamental representation; the resulting principal pairs correspond
to rank $n$ vortices. We then can take the set $S$ to be the set of
endomorphisms $\chi$ of the vector bundle $V=P_K\times_K\CC^n$ such
that:
\begin{enumerate}
\item the characteristic polynomial $P(\chi_x)$ of $\chi$ acting on the
fibre over $x\in X$ does not depend on $x$,
\item $P(\chi_x)$ has at most two different roots, and both 
come from the set $\{ 0, \pm \imag\}$, 
\item for any root $\alpha$ of $P(\chi_x)$ the set $\Ker(\chi-\alpha\Id)$
is in fact a subbundle of $V$.
\end{enumerate}
\end{example}
\begin{example}
If we take $\V=S^2\CC^n$ with the obvious linear action of
$K=\U(n)$ then the resulting principal pairs are the so called
quadric bundles. Then we can take $S$ to be the set of endomorphisms
which induce holomorphic filtrations of length at most three (i.e.,
the sections $\chi\in\Omega^0({\ad} P_K)$ belonging to $S$ have at most
three different eigenvalues). This is proved in
\cite{GGM}.
\end{example}
In both examples, the reason why one can restrict to filtrations of
lengths two and three is the following. Let us fix a pair $(A,\Phi)$
and consider a $\ov{\partial}_A$-holomorphic filtration\footnotemark
$$\VVV=(0\subset V_1\subset\dots\subset V_k=V)$$
\footnotetext{of vector subbundles if $\dim X=1$, or subsheaves in general
if $\dim X>1$} induced by any section of $\ad P_K$. Any choice of 
real numbers $\lambda=(\lambda_1<\dots<\lambda_k)$ determines a 
negative subbundle $V_{\lambda}^-\subset V$. Moreover, there is a 
unique section $\chi_{\VVV,\lambda}$ in $\Omega^0(\ad P_K)$ inducing 
the filtration $\VVV$ and having $\imag\lambda$ as set of 
eigenvalues. With 
$$f(\lambda)=\deg(\chi_{\VVV,\lambda})-
\int\la\chi_{\VVV,\lambda},c\ra\ ,$$
the stability condition becomes:
$f(\lambda)>0$ whenever $\Phi\subset V_{\lambda}^-$.
Now let
$$\Lambda=\Lambda(\VVV,\Phi)=\{\lambda\mid\Phi\subset
V_{\lambda}^-\}\ .$$ This is a convex polytope. Since $f(\lambda)$ is
a linear function, it is thus enough to check the stability condition
for the choices of $\lambda$ belonging to any subset
$\Lambda'\subset\Lambda$ whose convex hull is
$\Lambda$. To prove that SSC applies in the above two examples one has to
study
what are the possible sets $\Lambda$ which can appear and find in
each case suitably simple subsets $\Lambda'$. See the proof of
Proposition
\ref{prop:fixedE1HK} for more details.
\subsection{SSC in the $\HHH$ case} \label{subs:HSSC}
If SSC holds for $(K,P_K,\rho,\Y)$ and the following refinement of
(1) is true:
(1') any $\sigma\in\hlie$ can be written as $\sigma=s_1+\dots+s_r$,
where each $s_i$ belongs to $S\cap\hlie$,
\noi then we say that SSC holds for the action of $\HHH$. The following
generalization of Lemma \ref{lemma:SSC} holds then true.
\begin{lemma}
Let $(A,\Phi)\in\XXX(P_K,\YY)$ be any pair, and let
$c_{\HHH}\in\hlie$ be a central element. Then $(A,\Phi)$ is
$(\HHH,c_{\HHH})$-stable if and only if for any $X_0\subset X$ whose
complement has complex codimension at least two and any $s\in
S|_{X_0}\cap \hlie$ inducing a
$A$-holomorphic filtration we have
$$\deg(s)+\int_{X_0}\lambda(\Phi;s)-\int_{X_0}\la s,c_{\HHH}\ra >0.$$
\label{lemma:SSCH}
\end{lemma}
\section{Subgroups determined by the structure group}\label{sect:StrSub}
In many examples, the subgroups of the gauge groups come from
subgroups of the structure groups of the principal bundles.  In this
section we describe three types of such subgroups.  While there are
undoubtedly other mechanisms which can produce subgroups of this
kind, the three we describe account for a surprisingly broad range of
examples.
\subsection{ Reduction of the structure group}\label{sect:reduction}
Suppose that the principal $K$-bundle $P_K$ admits a reduction of the
structure group to $H$. Let $P_H$ be the corresponding principal
bundle. Thus $P_K=P_H\times_HK$ and the Adjoint bundle of
$P_H$, i.e.
\begin{equation}
{\Ad}P_H=P_H\times_{{\Ad}}H\ ,
\end{equation}
\noi is a subbundle of ${\Ad}P_K=P_K\times_{{\Ad}}K$.
The gauge group of $P_H$, i.e. $\GGG(P_H)=\Omega^0({\Ad}P_H)$, can thus
be viewed as a subgroup of $\GGG=\Omega^0({\Ad}P_K)$. The subgroup of
$\GGG$ corresponding to $H$ is given , in this case, by
\begin{equation}
\HHH=\GGG(P_H)
\end{equation}
\noi Similarly, the complexification $\HHC$\ is the (complex) gauge group
of the $H^{\CC}$-principal bundle $P_H\times_H H^{\CC}$.
\begin{lemma}  The subgroup $\HHH$ and its complexification
satisfy the conditions in section (\ref{subs:subgp}). If $A$\ is a
connection on $P_K$\ which comes from a connection on $P_H$, i.e. if
$A$\ is an $H$-connection, then
$d_A(\hlie)\subset
\hlie^1$.
\end{lemma}
\begin{pf} The first statement is obvious.
With respect to a local frame, any connection on $P_K$ can be written
as
$d_A=d+a$\ where $a$\ is a
1-form with values in the Lie algebra of $K$. If $A$\ is an
$H$-connection, then we can assume that $a$ takes its values in the Lie
algebra of
$H$. The result follows from this.
\end{pf}

Set
\begin{equation}
\XXX_{\HHH}=\{(A,\Phi)\in \XXX(P_K,\YY)\ |
A \mathrm{\ is\ an\ } H\mathrm{-connection}\}\ ,
\end{equation}
\begin{lemma} The subspace $\XXX_{\HHH}\subset\XXX(P_K,\YY)$ is invariant
under the action of $\HHH$ and also under the extension of this
action to $\HHC$
\end{lemma}
\begin{pf}  Since $A$ is an $H$-connection, it can be viewed as a
connection on $P_H$.   Moreover, $\YY_K$ can be described as a an
associated bundle to $P_H$, i.e. as $\YY=P_H\times_H \Y$.  Thus the
pair
$(A,\Psi)$ can equally well be treated as principal pairs of type
$(P_H,\YY)$.
The results follows from this.
\end{pf}

\noi Given a principal pair type $(P_K,\YY_K)$ in which the
structure group of $P_K$ reduces to $H$, we can thus apply the Main
Theorem to any simple pair $(A,\Psi)$\ in $\XXX_{\HHH}$. Such pairs
can, however, be viewed as a pair of type $(P_H,\YY_H)$. The two
points of view are equivalent:
\begin{prop}\label{prop:equiv} Treating the pair $(A,\Phi)$\ as a pair
of type $(P_H,\YY_H)$ is equivalent to treating it as an $\HHH$-pair
of type
$(P_K,\YY_K)$. More precisely: the stability notions coincide, as do the
generalized vortex equations, and the Hitchin--Kobayashi
correspondence for pairs of type $(P_H,\YY_H)$\ is equivalent to the
Main Theorem applied to $H$-pairs of type $(P_K,\YY_K)$.
\end{prop}
\begin{remark} Consider, for example, the case in which $H$ is a subgroup
of $\U(n)$ and $\Y=\CC^n$. The bundle $E=P_H\times_H\Y$ can then be
viewed as a rank $n$ vector bundle with a reduction of its structure
group from $\U(n)$ to $H$. In this case Proposition \ref{prop:equiv}
says that principal pairs of type
$(P_H,\YY_H)$ can equally well be viewed as special holomorphic pair on
the vector bundle $E=P_H\times\CC^n$; namely as holomorphic pairs
which are compatible with the reduction of structure group of $E$
from $\U(n)$ to $H$.
\end{remark}
\subsection{Normal subgroups}\label{sect:normal}
Let $H$  be a normal subgroup of $K$ and suppose that its
complexification
$H^{\CC}$ is a normal subgroup of $G$. Even if the
structure group of $P_K$ does not reduce to $H$ , we can define the
subbundle
\begin{equation}
P_K\times_{{\Ad}}H\subset {\Ad}P_K \ .
\end{equation}
We thus get a  subgroup
\begin{equation}
\HHH=\Omega^0(P_K\times_{{\Ad}}H)\subset \GGG\
\end{equation}
\noi with complexification
\begin{equation}
\HHC=\Omega^0(P_G\times_{{\Ad}}H^{\CC})\subset\GGG^{\CC}\ .
\end{equation}
\begin{lemma}
Subgroups of this sort satisfy the conditions (1)-(3) in section
(\ref{subs:subgp}). Furthermore, given any connection, $A$, on $P_K$,
we get $d_A(\hlie)\subset \hlie^1$.
\end{lemma}
\begin{pf}
The fact that $d_A(\hlie)\subset \hlie^1$ follows directly from the
fact that if $H\subset K$ is a normal subgroup, then the Lie algebra
of $H$ is an ideal in the Lie algebra of $K$.
\end{pf}

\begin{remark}  The main examples we have in mind (see Section
\ref{sect:examples}) all have real structure groups of the form
\begin{equation}
K=K_1\times K_2\times\dots\times K_n\ .
\end{equation}
The normal subgroups of $K$ are obtained by restricting to the identity
element in some of the factors.
\end{remark}
\subsection{Constant gauge transformations}\label{sectn:constant}
Suppose that $P_K=X\times K$, i.e. suppose that $P_K$ is the trivial
principal $K$-bundle. Then $K$ (or indeed any subgroup $H\subset K$)
embeds in $\GGG$ as the group of constant gauge transformations. In
this case, given $H\subseteq K$, we may take
\begin{equation}
\HHH =  H\ ,\ \HHC= H^{\CC}\ .
\end{equation}
\noi The requirements of section (\ref{subs:subgp}) are
 certainly satisfied by subgroups of this kind. We can define an
$\HHC$-invariant subspace in
$\XXX(P_K,\YY)$ by fixing the connection to be the trivial connection
on the trivial bundle, i.e. we can define
\begin{equation}
\XXX_{\HHH}=\{(0,\Phi)\in \XXX(P_K,\YY)\ |
0\ \mathrm{denotes\ the\ trivial\ connection}\}\ .
\end{equation}
\noi Since the connection in $\XXX_{\HHH}$ is trivial, the requirements of
Definition
\ref{defn:subgp} are satisfied by $(\HHH,\HHC,\XXX_{\HHH})$,
i.e. $d_A(\hlie)\subset \hlie^1$ for all connections which occur in
$\XXX_{\HHH}$.
\begin{remark} We can generalize this situation to the case in which
$K=K_1\times K_2$ (or indeed a product of more than two factors), and
$P_K$ is a fibre product (say $P_{K_1}\times
P_{K_2}$) with one of the factors trivial. An example of this sort
arises in the description of Coherent Systems (see Section
\ref{subs:CS}).
\end{remark}
\section{Examples}\label{sect:examples}\hfil
Subgroups of the gauge group occur naturally when the gauge group is
a product of two or more groups. In this section we describe some
examples of this sort. Thus we consider principal pairs for which the
complex structure group is a product, say
\begin{equation}\label{eqtn:Gproduct}
G=G_1\times G_2\times\dots\times G_p\ ,
\end{equation}
\noi and the principal $G$ bundle is a fibre product, say
\begin{equation}
\label{eqtn:PGproduct}
P_G=P_{G_1}\times P_{G_2}\times\dots\times P_{G_p}\ ,
\end{equation}
\noi where $P_{G_i}$ is a principal $G_i$-bundle.  If the compact real form
of $G_i$ is $K_i$, then the compact real
form of $G$ is
$$ K=K_1\times K_2\times\dots\times K_p\ .$$
\noi The real principal bundle obtained by a reduction of structure group
from $G$ to $K$ is thus
$$P_K=P_{K_1}\times P_{K_2}\times\dots\times P_{K_p}\ .$$
The natural subgroups of the gauge groups in such examples are
determined by restricting to the identity element in one or more of
the factors in $G_1\times G_2\times\dots\times G_p$.
To complete the specification of a principal pair type we need a
Kaehler manifold together with a Hamiltonian $K$-action. In this paper
we restrict ourselves to examples in which the Kaehler manifold is a
vector space $\V$.  In most of our examples $\V$ has the form
$$\V=\V_1\otimes \V_2\otimes\dots\otimes \V_p$$
\noi where for $1\le i\le p$ each $\V_i$ is a complex vector space of
dimension
$n_i$. Let $\la\ ,\ \ra_i$ be a hermitian inner product on $\V_i$ and
let $\omega_i$ be the corresponding Kaehler form. Thus
\begin{equation}
\omega_i(x,y)=\frac{1}{2\imag}(\la x,y\ra_i-\la y,x\ra_i)\ .
\end{equation}
\noi We let $\la\ ,\ \ra$ be the hermitian inner product on $\V$ determined
by
the inner products on the
$\V_i$. Thus if, for $1\le i\le p$, the collections
$\{e^{i}_1,e^{i}_2,\dots,e^{i}_{n_{i}}\}$ is a
unitary frame for $\V_i$, then the tensor products
$\{e^{1}_{i_1}\otimes e^{2}_{i_2}\otimes\dots\otimes e^{p}_{i_p}\}$ form
a unitary frame for $\V$. We let $\Omega$ be the corresponding Kaehler
form.
The principal pairs we consider are then of the form
$(P_K,\VV_K)$, with $\VV_K=P_K\times_{\rho}\V$. Furthermore in this section
we
consider only examples in which the $K$ action on $\V$ arises from  a
representations $\rho:K\longrightarrow \U$, where
$U$ denotes the group of unitary transformations on $(\V,\la\ ,\ \ra)$.
If we let $E_i=P_{K_i}\times_{\U_i}\V_i$ denote the vector bundle
associated to $P_{K_i}$ by the standard representation of $\U_i$ on
$\V_i$, then from the holomorphic point of view, the pairs on $(P_K,\VV_K)$
are equivalent to
\begin{itemize}
\item a collection of holomorphic bundles $\EEE_1,\dots,\EEE_p$ together
with
\item a section of an associated holomorphic bundle $\VVV_K$.
\end{itemize}
\noi Here $\EEE_i$ denotes the holomorphic bundle obtained by putting a
holomorphic structure on the smooth bundle $E_i$, and similarly
$\VVV_K$ denotes the holomorphic bundle obtained by putting a
holomorphic structure on the smooth bundle $\VV_K$.
The form of the generalized vortex equations on pairs of this type,
and in particular the equations which result from passing to a
subgroup of the gauge group (as in our Main Theorem), is the result
of the following observations. Notice that a connection on $P_K$ is
equivalent to a $p$-tuple of connections on $P_{K_1},\dots,P_{K_p}$.
Writing $A=(A_1,\dots,A_p)$, we see that the corresponding curvature
term in the generalized vortex equations (cf. (\ref{eqn:defxi})) is
of the form
$$\Lambda F_A=(\Lambda F_{A_1},\dots,\Lambda F_{A_p})\ ,$$
\noi where $\Lambda F_A$ takes its values in
$\Lie\GGG=\bigoplus\Lie\GGG_i$
(and each $\Lambda F_{A_i}$ has its values in $\Lie\GGG_i$).
The other term in the vortex equation is described in the following
lemma.
\begin{lemma}
\label{lemma:mmap7}
Let $\GGG_i$ be the gauge group for $P_{K_i}$, and
write $\GGG=\GGG_1\times\dots\times\GGG_p$.  Fix faithful unitary
representations $\rho_{a,i}:K\to \U(\WW_{a,i})$, where
$\WW_{a,i}$  are finite dimensional Hermitian vector spaces.  As in Section
\ref{subs:geomsetting} use these to define inner products on $\Lie \GGG_i$
and hence to get inclusions
$\Lie\GGG_i\subset(\Lie\GGG_i)^*$.  Let
$$\mu_i:\Omega^0(V_K)\longmapsto\Lie\GGG_i\subset(\Lie\GGG_i)^*$$
\noindent be the moment map for the action of $\GGG_i$.
Then the moment map for $\GGG$ is the map
$$\mu:\Omega^0(V_K)\longmapsto
\bigoplus_{i=1}^p\Lie\GGG_i\subset
\bigoplus_{i=1}^p(\Lie\GGG_i)^*=(\Lie\GGG)^*$$ given by
\begin{equation}
\mu(\Phi)=(\mu_1(\Phi),\dots,\mu_p(\Phi)).
\end{equation}
\end{lemma}
\subsection{Example 1 (Tensor product bundles)}\label{subs:tensor}
Consider the case\footnotemark
\footnotetext{There is an obvious generalization of this example to tensor
products with more than two factors. Since no new ideas are involved,
we discuss only the simplest case} where
\begin{itemize}
\item $G_i=\GL(n_i)$ and $K_i=\U(n_i)$ for $i=1, 2$, so
\item $G=\GL(n_1)\times \GL(n_2)$, and
correspondingly $K=\U(n_1)\times \U(n_2)$,
\item $\V_i=\CC^{n_i}$ for $i=1$ and $2$ and $\rho_i:G_i\longrightarrow
\GL(n_i)$
is the standard representation on $\V_i$.
\end{itemize}
We take $ \V=\V_1\otimes \V_2$ and let
$\rho=\rho_1\otimes\rho_2$  be the tensor product
representation of $G=\GL(n_1)\times \GL(n_2)$ on $\V$. Thus
$\VV_G=P_G\times_{\rho}\V$ is a tensor product of vector bundles,
i.e. $$\VV_G=\VV_1\otimes\VV_2\ ,$$
\noi where $\VV_i=P_{\GL(n_i)}\times_{\rho_i}\V_i$
is the rank $n_i$ vector bundle associated to $P_{\GL(n_i)}$. A
holomorphic structure on
$P_G$ is thus equivalent to a holomorphic structure on $\VV_G$
such that the resulting holomorphic vector bundle, $\VVV$, is a
tensor product of holomorphic bundles, i.e.
$\VVV=\VVV_1\otimes \VVV_2$ (where $\VVV_i$ denotes
holomorphic bundle obtained by putting a holomorphic structure on
$\VV_i$). Hence
{\it With $P_G$ and $\VV_G$ as above, a principal pair of type $(P_G,
\VV_G)$ is equivalent to a pair
$(\VVV,\Phi)$, where $\VVV$ is holomorphic bundle of the form
$\VVV_1\otimes \VVV_2$ and $\Phi$ is a holomorphic section of $\VVV$.}
We now show how our Main Theorem leads to a Hitchin--Kobayashi
correspondence for pairs on
$\VVV_1\otimes\VVV_2$ but for which one of the tensor factors is
regarded as fixed.\footnotemark
\footnotetext{Of course we can identify $\VVV_1\otimes\VVV_2$ with
$\Hom(\VVV_2^*,\VVV_1)$ (or, for that matter, with
$\Hom(\VVV_1^*,\VVV_2$). This may tempt one to interpret a pair
$(\VVV_1\otimes\VVV_2,\Psi)$
as a holomorphic triple (as in \cite{BGPtriples}) on
$(\VVV_1,\VVV_2^*)$, but one should resist the temptation. The two
types of augmented bundle differ precisely in the fact that the
bundles underlying the pair
$(\VVV_1\otimes\VVV_2,\Psi)$ are $\VVV_1$ and $\VVV_2$, while those
underlying the triple are $\VVV_1$ and $\VVV^*_2$. We will return to
triples in the next section.}
The complex gauge group $G=\GL(n_1)\times \GL(n_2)$ has a normal
subgroup defined by
\begin{equation}
H^{\CC}=\GL(n_1)\times \{1\} \ .
\end{equation}
\noi Let $\HHC$ be the corresponding subgroup of
$\GGC$, as in Section \ref{sect:normal}. This is the complexification
of the subgroup $\HHH\subset\GGG$ defined by
$H=\U(n_1)\times\{1\}$.
The corresponding Lie algebras are the following:
$$\Lie\GGG = \{(s_1,s_2)\in\Omega^0(\End V_1)\oplus\Omega^0(\End V_2)
\mid s_1+s_1^*=0,\ s_2+s_2^*=0\ \},$$
and
$$\hlie = \{(s_1,s_2)\in\Lie\GGG\mid s_2=0\}.$$
Then we take as projection the map $\pi_{\hlie}:\Lie\GGG\to\hlie$
which sends $(s_1,s_2)$ to $(s_1,0)$.
If we fix a connection, say $A_2^0$, on
$P_{\U(n_2)}$, then we get a subspace
$\XXX_H\subset\XXX(P_K,E_K)$ defined by the condition that
the connections on $P_K$ are of the form $A=(A_1 , A_2^0)$. (As
connections on the vector bundle $V=V_1\otimes V_2$, these correspond
to connections of the form  $A=A_1\otimes I_2 + I_1\otimes A_2^0$.)
\begin{lemma}  The subspace $\XXX_H$ is $\HHH$-invariant, and the
corresponding complex subspace of
$\XXX(P_G,E_G)$ is $\HHC$-invariant. The data $(\HHH,\HHC,\XXX_H)$
 determines  subgroup settings.
\end{lemma}
Since from the holomorphic point of view the bundle $\VVV_2$ factor
in $\VVV_1\otimes\VVV_2$ is fixed (namely, it is the bundle obtained
by considering on $V_2$ the $\ov{\partial}$-operator given by
$A_2^0$), we say that the subspace
$\XXX_{H}$ defines the configuration space of
\it $\VVV_2$-twisted pairs on $V_1$\rm. The
$\HHC$ orbits correspond to the isomorphism classes of such twisted
pairs.
\begin{prop}
\label{prop:fixedE1HK}
Let $\HHH$ and $\XXX_H\subset\XXX(P_K,V_K)$ be as above. Define the
auxiliary representation
$\rho_a:K\longrightarrow \U(\CC^{n_1}\oplus\CC^{n_2})$
using the standard representations, and thereby fix an inclusion
$\Lie\GGG\subset(\Lie\GGG)^*$. Let
$c_{\HHH}:=-\sqrt{-1}(cI_1, 0)\in\hlie$,
where $c\in\RR$ is any number and $I_1\in\End V_1$ is the identity. Then
\begin{enumerate}
\item Fix a system of local unitary frames for $\VV_2$, say
$\{e^{\alpha}_i\}_{i=1}^{n_2}$ over $U_{\alpha}\subset X$, where
$\{U_{\alpha}\}$ is a suitable open cover of $X$. For any $x\in U_{\alpha}$
write
$$\Phi(x)=\sum_{i=1}^{n_2}\varphi^{\alpha}_i(x)\otimes
e^{\alpha}_i(x)$$
\noi where the $\varphi^{\alpha}_i(x)$ are locally defined
sections of $\VV_1$.
When applied to any $((A_1 , A_2^0), \Phi)\in\XXX_H$, and for any
$x\in U_{\alpha}$, the
$(\HHH,c_{\HHH})$-vortex equation takes the form
\begin{equation}
\label{eqtn:H1c}
\sqrt{-1}\Lambda F_{A_1}(x) +
\sum_{i=1}^{n_2}\varphi^{\alpha}_i(x)\otimes (\varphi^{\alpha}_i(x))^*=
cI_1.
\end{equation}
\item Let $((A_1 , A_2^0),\Phi)\in\XXX_H$ and denote by
$\VVV_1$ the holomorphic bundle obtained by considering on $V_1$
the $\ov{\partial}$-operator given by $A_1$. Assume that $\Vol X=1$
(the general case can be reduced to it by rescaling appropriately).
\begin{enumerate}
\item The pair $((A_1 , A_2^0),\Phi)$ is
$(\HHH,c_{\HHH})$-stable if and only if:
all coherent subsheaves of $\VVV_1\otimes\VVV_2$
of the form $\VVV'\otimes\VVV_2$ satisfy:
$$\mu(\VVV')<c$$
(where as usual \footnotemark\footnotetext{We apologize for the 
excessive use of the letter $\mu$}, 
 $\mu(\VVV)=\deg(\VVV)/\rk(\VVV)$),
\noi and if $\Phi\in H^0(\VVV'\otimes\VVV_2)$ then
$$\mu(\VVV_1/\VVV')>c.$$
\item The pair $((A_1 , A_2^0),\Phi)$ is simple if and only if
there is no
holomorphic nontrivial splitting $\VVV_1=\VVV_1'\oplus\VVV_2''$
such that $\Phi\in H^0(\VVV_1'\otimes\VVV_2)$.
\end{enumerate}
\end{enumerate}
\end{prop}
\begin{remark}
Though the local sections $\varphi^{\alpha}_i$ depend on
the choice of local unitary frames $\{e^{\alpha}_i\}_{i=1}^{n_2}$,
the expressions $\sum_{i=1}^{n_2}\varphi^{\alpha}_i(x)\otimes
(\varphi^{\alpha}_i(x))^*$ do not and are thus globally defined. This
can be checked directly as follows. Let $
\{\tilde{e}^{\alpha}_i\}_{i=1}^{n_2}$ be another unitary frame for $\VV_2$,
related to
$\{e^{\alpha}_i\}_{i=1}^{n_2}$ by
$\tilde{e}^{\alpha}_i=T_{ji}e^{\alpha}_j$.
Since both frames are unitary, the elements $T_{ji}$
define a unitary matrix. Writing
$\Phi(x)=\sum_{i=1}^{n_2}\tilde{\varphi}^{\alpha}_i(x)\otimes
\tilde{e}^{\alpha}_i(x)=\sum_{i=1}^{n_2}T_{ji}\tilde{\varphi}^{\alpha}_i(x)
\otimes e^{\alpha}_i(x)$, we
see that $\tilde{\varphi}^{\alpha}_i=T{ji}\varphi^{\alpha}_j$. Thus
\begin{equation}
\sum_{i=1}^{n_2}\varphi^{\alpha}_i(x)\otimes
(\varphi^{\alpha}_i(x))^*= \sum_{j,k=1}^{n_2}\sum_{i=1}^{n_2}
T_{ji}{T^*_{ik}}\tilde{\varphi}^{\alpha}_j(x)\otimes
(\tilde{\varphi}^{\alpha}_k(x))^*\ =
\sum_{j=1}^{n_2}\tilde{\varphi}^{\alpha}_j(x)\otimes
(\tilde{\varphi}^{\alpha}_j(x))^*
\end{equation}
\noi where we have used the unitarity of $T$ in the last equality.
\end{remark}
\begin{pf} {\it (of Proposition \ref{prop:fixedE1HK})}
We first analyze the equations. Recall from (\ref{eqtn:HcVortex})
that the $(\HHH,c_{\HHH})$ vortex equations for a point 
$(A,\Phi)$ in $\XXX_H$ are given by
\begin{equation}\label{eqtn:HCvortex2}
\pi_{\hlie}(\Lambda F_{A}+\mu(\Phi))=c_{\HHH}\ .
\end{equation}
\noi Given a connection $A=(A_1,A_2)$ on
$P_G=P_{\GL(n_1)}\times P_{\GL(n_2)}$, we get
\begin{equation}
\Lambda F_A=(\Lambda F_{A_1},\Lambda F_{A_2}),
\end{equation}
\noi hence $\pi_{\hlie}(\Lambda F_A)=\Lambda F_{A_1}$.
The term $\mu(\Phi)$ in (\ref{eqtn:HCvortex2})
is determined by the $\U(n_1)\times \U(n_2)$ moment map on
$\CC^{n_2}\otimes\CC^{n_1}$, with respect to the usual Kaehler
structure. Denoting this too by $\mu$, we have
\begin{equation}
\mu(p)=(\mu_1(p),\mu_2(p))
\end{equation}
\noi where for $i=1,2$, the maps
$\mu_i:\CC^{n_2}\otimes\CC^{n_1}\longmapsto \mathfrak{u}(n_i)^*\cong 
\mathfrak{u}(n_i)$ are the
moment maps for the action of $\U(n_1)=\U(n_1)\times \{1\}$ and
$\U(n_2)=  \{1\}\times \U(n_1)$ respectively.
Thus, as in Lemma \ref{lemma:mmap7},
\begin{equation}
\mu(\Phi)=(\mu_1(\Phi),\mu_2(\Phi))\ .
\end{equation}
Consequently, $\pi_{\hlie}(\mu(\Phi))=\mu_1(\Phi).$ To obtain
(\ref{eqtn:H1c}) it thus remains to evaluate the $\U(n_1)$ moment map
$\mu_1$. By Propositions \ref{prop:A1} and \ref{prop:muVK}
in the Appendix, this can be described as
follows. Fixing a unitary basis, say $\{e_i\}_{i=1}^{n_2}$, for
$\CC^{n_2}$, we can write any
$x\in\CC^{n_1}\otimes\CC^{n_2}$ as $x=\sum_{i=1}^{n_2}x_i\otimes e_i$, where
the $x_i$ are vectors in $\CC^{n_1}$. The moment map $\mu_1(x)$ is
given by
 \begin{equation}
 \label{eqtn:U1mmap}
 \mu_1(\sum_{i=1}^{n_2}x_i\otimes e_i)=-\sqrt{-1}\sum_{i=1}^{n_2}x_i\otimes
 x^*_i\ ,
 \end{equation}

 \noi where $x^*_i$ denotes the dual element in $(\CC^{n_1})^*$.
 Equivalently, representing vectors in $\CC^{n_1}$ by column
 vectors and using row vectors to represent their duals,  the
 terms $x_i\otimes x^*_i$ become matrices $x_i\overline{x_i}^t$ in
 $\Hom(\CC^{n_1},\CC^{n_2})$.
Notice that the formula in (\ref{eqtn:U1mmap}) is independent of the
choice of unitary basis for $\CC^{n_2}$, and is also $\U(n_1)$
invariant. In order to describe the corresponding map
\begin{equation}
\mu_1:\Omega^0(\VVV_1\otimes\VVV_2)\longmapsto \Lie\GGG_1
\end{equation}
\noi we may thus pick local frames for the bundles and apply
(\ref{eqtn:U1mmap}) directly to these. This leads immediately to the
term $\sum_{i=1}^{n_2}\varphi^{\alpha}_i(x)\otimes
(\varphi^{\alpha}_i(x))^*$ in (\ref{eqtn:H1c}).

\par We now show what stability with respect to the subgroup
$\HHH$ means. For convenience we will assume that $\dim X=1$
(in the general case we should take into account that the filtrations
appearing in the definition of stability could be in principle defined
only in the complementary of analytic subsets of $X$ of codimension
$\geq 2$, hence we should consider  reflexive subsheaves and not
only about subbundles).
By Definition \ref{defn:HcStability} (together with \S
\ref{sss:vectorbdl}) the pair
$((A_1,A_2^0),\Phi)$ is $(\HHH,c_{\HHH})$-stable if and only if, for
any $\chi\in\hlie$ inducing a holomorphic filtration such
that $\Phi\subset (\VVV_1\otimes \VVV_2)^-(\chi)$, the following
inequality holds
\begin{equation}\label{eqtn:chige0}
\deg\chi-\int_X\la \chi,c_{\HHH}\ra>0\ ,
\end{equation}
where $\deg\chi$ is given in Definition \ref{defn:degchi}. This can
be re-formulated more concretely as follows. We observe that any
$\chi\in\hlie$ is of the form $(\chi_1, 0)$. If $\chi$ induces a
$A$-holomorphic filtration of
$W=P_K\times_{\rho_a}\U(\CC^{n_1}\oplus\CC^{n_2})$
then it has fibrewise constant eigenvalues, among which zero is
included. Suppose the eigenvalues of $\-\sqrt{-1}\chi$ are the real
numbers $\alpha_1\le\alpha_2\le\dots\le\alpha_r$\footnotemark, with
$\alpha_j=0$. If the filtration of $W$ is
$$0\subset W_0\subset W_1\subset\dots\subset W_r$$
then
\begin{equation}
W_i= \left\{
\begin{matrix}
\VVV_{1,i}\oplus 0& \ \mathrm{if}\ 1\le i\le j-1\\
\VVV_{1,i}\oplus \VVV_2& \ \mathrm{if}\ i\ge j
\end{matrix}\right\}
\end{equation}
\noi where the $\VVV_{1,i}$ are the terms in the
holomorphic filtration
\begin{equation}
0\subset \VVV_{1,1}\subset\dots\subset \VVV_{1,r}=\VVV_1\ .
\end{equation}
\footnotetext{the fact that the inequalities are
not strict means that one can obtain some $\chi$ with different
choices of filtrations and weights}
\noi determined by $\chi_{1}\in\Lie\GGG_1$.
We deduce that
\begin{align*}
\deg\chi=& \alpha_r\deg W
+\sum_{k=1}^{r-1}(\alpha_k-\alpha_{k+1})\deg W_k\\ =& \alpha_r(\deg
\VVV_1+\deg\VVV_2)
+\sum_{k=1}^{j-1}(\alpha_k-\alpha_{k+1})\deg\VVV_{1,k}
+\sum_{k=j}^{r-1}(\alpha_k-\alpha_{k+1})(\deg\VVV_{1,k}+\deg\VVV_2)\\
=& \alpha_r\deg \VVV_1
+\sum_{k=1}^{r-1}(\alpha_k-\alpha_{k+1})\deg\VVV_{1,k}+\alpha_j\deg\VVV_2\\
=& \alpha_r\deg \VVV_1
+\sum_{k=1}^{r-1}(\alpha_k-\alpha_{k+1})\deg\VVV_{1,k}
\end{align*}
\noi where in the last line we have used the fact that $\alpha_j=0$.
Also, since $\Vol X=1$, we get
\begin{equation}
\int_X\la\chi,c_{\HHH}\ra=\int_X\la\chi_1,c\ra=
c(\alpha_r\rk \VVV_1+\sum_{k=1}^{r-1}(\alpha_k-\alpha_{k+1})
\rk\VVV_{1,k})\ .
\end{equation}
Thus the condition (\ref{eqtn:chige0}) is equivalent to
\begin{equation}
 \label{eqtn:deg-alpha-0}
\deg(\alpha)>0\
\end{equation}
\noi where
\begin{equation}
\label{eqtn:deg-alpha-defn}
\deg(\alpha):=\alpha_r(\deg \VVV_1-c\rk\VVV_1)+
\sum_{k=1}^{r-1}(\alpha_k-\alpha_{k+1})
(\deg\VVV_{1,k}-c\rk\VVV_{1,k})\ .
\end{equation}

%
To determine which filtrations must satisfy (\ref{eqtn:deg-alpha-0})
we need to identify when the condition
$\Phi\in H^0((\VVV_1\otimes\VVV_2)^-(\chi))$ is satisfied.
Since the action of
$\rho(\chi)$ on $\VVV_1\otimes\VVV_2$ is given by
$\chi_{1}\otimes\Id_{\VVV_2}$, its eigenvalues are
$\{-\imag\alpha_1,\dots,-\imag\alpha_r\}$,
occurring with multiplicities determined by the
ranks of  $\{\VVV_{1,1}\otimes\VVV_2,
\VVV_{1,2}\otimes\VVV_2,\dots, \VVV_{1,r}\otimes\VVV_2\}$.

In this case \footnotemark
\footnotetext{In later examples the corresponding integers will be
defined slightly differently} we thus associate to $\chi$ the
integers
\begin{align}\label{eqn:palpha}
p(\alpha)& =\max\{i\mid\alpha_i\leq 0\ \}\\ p(\chi)& =
\min\{i\mid\Phi\ \in  H^0(\VVV_{1,i}\otimes\VVV_2)\ \}
\end{align}
\noi Then $(\VVV_1\otimes\VVV_2)^-(\chi)= \VVV_{1,p(\alpha)}
\otimes\VVV_2$ and furthermore, $\Phi \in
H^0((\VVV_1\otimes\VVV_2)^-(\chi))$ if and only if
$p(\chi)\le p(\alpha)$. The stability condition can thus be
re-formulated as follows:

\noi {\it A pair corresponding to a holomorphic bundle $\VVV_1\otimes\VVV_2$
and a section $\Phi\in H^0(\VVV_1\otimes\VVV_2)$ is
$(\HHH,c_{\HHH})$-stable if and only if the following holds:  Take any
holomorphic filtration
$0\subset \VVV_{1,1}\subset\dots\subset \VVV_{1,r}=\VVV_1$.  Set
$$\Lambda=\{\alpha\in\RR^r\mid \alpha_i\leq\alpha_{i+1}\text{ for
$1\leq i\leq r-1$ and } p(\alpha)\ge p(\chi) \}\ ,$$
\noi where $p(\chi)$ and $p(\alpha)$ are as above. Then,
for any $\alpha\in\Lambda$ we have
$$\deg(\alpha)>0\ $$
\noi where $\deg(\alpha)$ is as in (\ref{eqtn:deg-alpha-defn})}.

We now show that the stability simplification condition applies.
Given a holomorphic filtration $\chi=(\chi_1,0)$, let $p=p(\chi)$.
Let $e_1,\dots,e_r$ be the canonical basis of $\RR^r$. Define, for
any $1\leq i\leq p$, $f_i=-\sum_{k\leq i} e_k$ and, for any
$p< j\leq r$, $g_j=\sum_{k\geq j} e_k$. It is straightforward to
check that $$\Lambda=\RR_{\geq 0}f_1+\dots+\RR_{\geq 0}f_r+
\RR_{\geq 0}g_{p+1}+\dots+\RR_{\geq 0}g_r.$$
Also, for any $\alpha\in\Lambda$, we  get
\begin{align}
\chi_{\alpha}& =\sum_{i=1}^p x_i\chi(f_i) +\sum_{i=p+1}^r y_i\chi(g_i)\\
\deg(\alpha)&  =\sum_{i=1}^p x_i\deg(f_i) +\sum_{i=p+1}^r y_i\deg(g_i).
\end{align}
\noi Here $\chi(f_i)=(\chi_1(f_i),0)$ where $\chi_1(f_i)$ denotes the
element in $\Lie \GGG_1$ whose eigenvalues are $-\sqrt{-1}\{-1, 0\}$,
with the multiplicity of the first being
$\sum_{k=1}^i\rank(\VVV_{1,k})$. Similarly
$\chi(g_i)=(\chi_1(g_i),0)$ where $\chi_1(g_i)$ denotes the element in $\Lie
\GGG_1$ whose eigenvalues are
$-\sqrt{-1}\{0,1\}$, with the multiplicity of the non-zero eigenvalue
being $\sum_{k=i}^r\rank(\VVV_{1,k})$.
\noi As in Example \ref{ex:SSC1} we define a subset $S\subset\hlie$
by the conditions
\begin{enumerate}
\item the characteristic polynomial $P((\chi_1)_x)$ of $\chi_1$ acting on
the
fibre over $x\in X$ does not depend on $x$,
\item $P((\chi_1)_x)$ has at most two different roots, both of which
come from the set $\imag\{-1, 0, 1\}$,
\item for any root $\alpha$ of $P((\chi_1)_x)$ the set
$\Ker(\chi_1-\alpha\Id)$
is in fact a subbundle of $\VVV_1$.
\end{enumerate}
\noi Then $\chi(f_i)$ and $\chi(g_i)$ are in $S$ and
the above computations show that condition (SSC1) in
Definition\ref{defn:SSC} applies. To verify the condition (SSC2') we
use the fact that the eigenvalues for $\chi(f_i)$ are both
non-positive while those for $\chi(g_i)$ are $0$ and $1$. Thus
$(\VVV_1\otimes\VVV_2)^-(\chi(f_i))=\VVV_1\otimes\VVV_2$, while
$(\VVV_1\otimes\VVV_2)^-(\chi(g_i))=\VVV_{1,i-1}\otimes\VVV_2$.
Notice in particular that
$(\VVV_1\otimes\VVV_2)^-(\chi(g_{p+1}))=(\VVV_1\otimes\VVV_2)^-(\chi)$.
It remains to interpret the stability condition as applied to
elements in $S$.  Consider an element $\sigma\in S$ which defines a
filtration $0\subset
\VVV'\subset\VVV_1$. Suppose first that the eigenvalues
are $\alpha_1=0$ and $\alpha_2=1$(i.e. an element of the form
$\chi(f_i)$). Then $p(\sigma)=1$ or $2$ (depending on whether $\Phi\in
H^0(\VVV'\otimes\VVV_2)$ or not) but $p(\alpha)=2$. Thus
$p(\sigma)$ is always less than or equal to $p(\alpha)$.
Moreover,
\begin{equation}\label{eqtn: fi-stable}
\deg(\sigma)>0 \Longleftrightarrow \mu(\VVV')<c\ .
\end{equation}
If $\alpha_1=-1$ and  $\alpha_2=0$, then $p(\sigma)\le p(\alpha)$ if and
only if $\Phi\in H^0(\VVV'\otimes\VVV_2)$, and
\begin{equation}\label{eqtn: gi-stable}
\deg(\sigma)>0 \Longleftrightarrow \mu(\VVV_1/\VVV')>c\ .
\end{equation}
The description of stability in (a) follows directly from these
observations. The characterization of simple pairs in (b) is
straightforward and is left as an exercise to the reader.
\end{pf}

\subsection{Example 2 (Fixed-$\EEE_2$ Triples)}\label{subs:triples}
In this example we take
\begin{itemize}
\item $G=\GL(n_1)\times \GL(n_2)$,
\item $K=\U(n_1)\times \U(n_2)$,
\item $\V_1=\CC^{n_1}$ and $\rho_1:\GL(n_1)\longrightarrow \GL(\V_1)$
 is the standard representation
\item $\V_2=\CC^{n_2}$, but $\rho_2:\GL(n_2)\longrightarrow \GL(\V_2)$
is the dual representation, i.e.
$$\rho_2(C)\cdot v=(C^{-1})^tv\ .$$
\end{itemize}
We take $ \V=\V_1\otimes \V_2$ and let
$\rho=\rho_1\otimes\rho_2$  be the tensor product
representation of $G=\GL(n_1)\times \GL(n_2)$ on $\V$. Equivalently, we
can take $\V=\V_1\otimes \V^*_2=\Hom(\V_2,\V_1)$, and regard $\rho$ as
the representation $\rho:G\longrightarrow \GL(\V)$ given by
\begin{equation}
\rho(C_1,C_2)(T)=C_1\circ T\circ C_2^{-1}\ .
\end{equation}
\begin{remark} If we set $E_i=P_{\GL(n_i)}\times_{\rho_i}\V_i$ for $i=1,2$,
i.e. if we let $E_i$ be the rank $n_i$ vector bundle associated to
$P_{\GL(n_i)}$ then with $P_G$ and $\V$ as above, {\it a principal pair of
type
$(P_G,\VV_G)$ is equivalent to the triple
$(\EEE_1, \EEE_2,\Phi)$, where $\EEE_i$ is holomorphic bundle with the
topological type of $E_i$, and $\Phi$ is a holomorphic section of
$\EEE_1\otimes\EEE_2^*=\Hom(\EEE_2,\EEE_1)$}. That is,  principal pairs of
this sort correspond to holomorphic triples in the sense of
\cite{BGPtriples}.
\end{remark}
Exactly as in \S \ref{sect:normal} we have two natural subgroups of
$K$, namely those corresponding to the normal subgroups
\begin{eqnarray}
H_1= \U(n_1)\times \{1\}\ , \\ H_2=\{1\}\times \U(n_2)
\end{eqnarray}
\noi in $K$. We denote the resulting subgroups of $\GGG$ by
$\HHH_1$ and
$\HHH_2$, and denote their Lie algebras by $\hlie_1$ and $\hlie_2$.
The complexifications, i.e. the subgroups of $\GGC$ corresponding to
the subgroups $\GL(n_1)\times\{1\}$ and $\{1\}\times \GL(n_2)$, are
denoted by  ${\HHH_1}^{\CC}$\ and ${\HHH_2}^{\CC}$,
with the obvious Lie algebras.
As in the previous example, if we fix an integrable connection, say
$A^0_2$, on $P_{\U(n_2)}$ then we can define a subspace
\begin{equation}
\label{eqn:XXX1}
\XXX_1=\{(A,\Phi)\in \XXX(P_K,\VV_K)\ |
A=(A_1, A^0_2)\}\
\end{equation}
\noi in $\XXX(P_K,\VV_K)$.  Equivalently, using the correspondence
between integrable connections on $P_K$ and holomorphic structures on
$P_G$, $\XXX_1$ can be regarded as a subspace of $\XXX(P_G,\VVV_G)$.
\begin{prop}\label{prop: fixed2}
The subspace $\XXX_1$ is an  $\HHH_1$-invariant subset in
$\XXX(P_K,\VV_K)$ and the corresponding subspace of $\XXX(P_G,\VV_G)$
 is an ${\HHH_1}^{\CC}$-invariant subspace.  The data
 $(\HHH_1,{\HHH_1}^{\CC},\XXX_1)$ determines a subgroup setting.
\end{prop}
\begin{remark}
If we let $P_{\GL(n_2)}$ have the holomorphic structure determined by
$A^0_2$ on $P_{\U(n_2)}$, and let $E_2$ have the holomorphic structure
determined by that on $P_{\GL(n_2)}$, then the
${\HHH_1}^{\CC}$-orbits in $\XXX_1$ correspond to isomorphism classes of
triples $(\EEE_1,\EEE_2,\Phi)$ in which the bundle $\EEE_2$ is fixed.
\end{remark}
\begin{remark}  There is, of course, an analog to Proposition
\ref{prop: fixed2} for fixed
$E_1$-triples, i.e. in which the roles of
$P_{\U(n_1)}$ and $P_{\U(n_2)}$ are reversed. In both cases, we get subgroup
settings
(as in Definition \ref{defn:subgp}), namely
$(\HHH_1,{\HHH_1}^{\CC},\XXX_1)$ for fixed $E_1$ triples and
$(\HHH_2,{\HHH_2}^{\CC},\XXX_2)$ for fixed $E_2$ triples. Our Main Theorem thus
applies. Formulated in terms of the triples data, it yields the
Hitchin--Kobayashi correspondence for fixed $E_1$ or fixed $E_2$
triples. We will state only the fixed $E_2$ case, the other being
entirely analogous.
\end{remark}
\begin{remark}  For obvious reasons triples such as those described in
Proposition
\ref{prop: fixed2}, i.e. in which one of the holomorphic bundles is
regarded as fixed, are called {\it fixed triple}. More specifically,
the objects described in that proposition may be called {\it fixed
$\EEE_2$-triples}. If the fixed bundle is
$\EEE_1=\FFF$, then the triple
$(\FFF,\EEE_2,\Phi)$ can be described as a holomorphic bundle (i.e.
$\EEE_2$) together with a morphism to the fixed bundle (i.e.
 $\Phi: \EEE_2\longrightarrow \FFF$). Since these constitute a special case
of
the framed modules studied by Huybrechts and Lehn \cite{HL}, we refer to them
as {\it framed bundles}. If we fix the bundle $\EEE_2$, we obtain an
object which can be described as a bundle together with a morphism
\it from\rm\ a fixed holomorphic bundle. Objects of this type, in which
the fixed bundle is a trivial rank $k$ bundle, provide a description
of coherent systems (see Section \ref{subs:CS}).
\end{remark}
\begin{prop}\label{thm:fixedE1HK.2}
 Let $\XXX_1\subset\XXX(P_K,\VV_K)$ be as above, and let
$\HHH:=\HHH_1$ be the subgroup corresponding to
$\U(n_1)\times\{1\}\subset K$. Let $\hlie=\Lie\HHH$ and define
$\pi_{\hlie}:\Lie\GGG\to\hlie$ to be the projection to the first
factor.
Use the standard representations to define an auxiliary representation
$\rho_a:K\longrightarrow \U(\CC^{n_1}\oplus\CC^{n_2})$, and thereby
get an inclusion $\Lie\GGG\subset(\Lie\GGG)^*$
(as in Section \ref{subs:geomsetting}).
Let $c_{\HHH}:=-\sqrt{-1}(cI_1, 0)$, where $c\in\RR$ is any real number
(then $c_{\HHH}$ is a central element in $\hlie$).
\begin{enumerate}
\item When applied to any $((A_1, A^0_2),\Phi)\in\XXX_1$, the
$(\HHH,c_{\HHH})$-vortex equation takes the form
\begin{equation}
\label{eqtn:HH1c}
\sqrt{-1}\Lambda F_{A_1}+\Phi\Phi^*= c I_1.
\end{equation}
\item Suppose that the volume $\Vol X=1$.
A point $((A_1, A^0_2),\Phi)\in\XXX_1$ is
$(\HHH,c_{\HHH})$-stable if and only if the corresponding
triple $(\EEE_1,\EEE_2,\Phi)$ satisfies the following condition: for
all coherent subsheaves $\EEE'_1\subset\EEE_1$
$$\mu(\EEE_1') < c\ ,$$
\noi and if $\Phi(\EEE_2)\subset\EEE'_1$ then
$$\mu(\EEE_1/\EEE'_1)>c.$$
The pair $((A_1, A^0_2),\Phi)\in\XXX_1$ is simple if there is no
holomorphic nontrivial splitting
$\EEE_1=\EEE_1'\oplus\EEE_1''$ such that $\Phi(\EEE_2)\subset\EEE_1'$.
\end{enumerate}
\end{prop}
\begin{pf} We first analyze the equations. The proof of (1) is essentially
the
same as the proof of (1) in Proposition  \ref{prop:fixedE1HK}, except
that in this case the term $\mu(\Phi)$ in  equation
(\ref{eqtn:HCvortex2}) is determined by the moment map for
$\U(n_1)\times \U(n_2)$ on $\Hom(\CC^{n_2},\CC^{n_1})$ (with respect to
the usual Kaehler structure). But, denoting this too by $\mu$, we have
\begin{equation}
\mu(T)=-\sqrt{-1}(T^*T,-TT^*)\ .
\end{equation}
\noi where $T$ is in $\Hom(\CC^{n_2},\CC^{n_1})$. The result now follows as
in
Proposition \ref{prop:fixedE1HK}.
To prove the statement on stability we follow exactly the same scheme
as in the proof of Proposition \ref{prop:fixedE1HK}. We consider
elements $\chi_{\alpha}=(\chi_{1,\alpha},0)\in\hlie$ with eigenvalues
$-\sqrt{-1}\{\alpha_1,\dots,\alpha_r\}$, and corresponding
filtrations $0\subset\EEE_{1,1}\subset\dots\subset\EEE_{1,r}$. We
identify filtrations $\chi(f_i)$ and $\chi(g_i)$, define the set
$S\subset \hlie$ in the same way as before, and verify  that SSC applies.
Notice that the eigenvalues for $\rho(\chi_{\alpha})$ on
$\Hom(\EEE_2,\EEE_1)$ are still
$-\sqrt{-1}\{\alpha_1,\dots,\alpha_r\}$. The filtration
determined by the corresponding eigen-subbundles is
$$0\subset \Hom(\EEE_2,\EEE_{1,1})\subset
\Hom(\EEE_2,\EEE_{1,2})\subset\dots\Hom(\EEE_2,\EEE_{1,r})
=\Hom(\EEE_2,\EEE_1)\ .$$ We thus define $p(\alpha)$ exactly as in
(\ref{eqn:palpha}), and set
\begin{equation}
p(\chi)=
\min\{i\mid\Phi\ \in  H^0(\Hom(\EEE_2,\EEE_{1,i}))\ \}
\end{equation}
\noi In particular,
\begin{itemize}
\item $p(\alpha(f_i))=2$ while $p(\chi(f_i))$ is either
one or two. Hence $p(\chi(f_i))\le p(\alpha(f_i))$ is always
satisfied.
\item $p(\alpha(g_i))=1$ while $p(\chi(g_i))$ is either
one or two. Hence $p(\chi(g_i))\le p(\alpha(g_i))$ is satisfied if 
$\Phi\in H^0(\Hom(\EEE_2,\EEE_{1,i}))$, i.e. 
$\Phi(\EEE_2)\subset\EEE_{1,i}$.
\end{itemize}
The statements in (2) of the Proposition now follow, as in the proof
of Proposition \ref{prop:fixedE1HK}, by evaluating the condition
$\deg(\sigma)>0$ for those elements $\sigma\in S$ which are of the
type defined by $\chi(f_i)$ and $\chi(g_i)$.
\end{pf}

\begin{remark} We can reformulate the stabiltiy condition in Proposition
\ref{thm:fixedE1HK.2}(2) as follows. For any $\alpha \in \mathbf{R}$
 define the
\emph{$\alpha$-degree} and
\emph{$\alpha$-slope} of a triple $(\EEE_1,\EEE_2,\Phi)$  to be
defined to be
\begin{align*}
  \deg_{\alpha}(\EEE_1,\EEE_2,\Phi)
  & = \deg(E_{1}) + \deg(E_{2}) + \alpha
  \rk(E_{2}), \\
  \mu_{\alpha}(\EEE_1,\EEE_2,\Phi)
  & =
  \frac{\deg_{\alpha}(\EEE_1,\EEE_2,\Phi)}
  {\rk(E_{1})+\rk(E_{2})}
  = \frac{\deg(E_{1} \oplus E_{2})}{\rk(E_{1})+
    \rk(E_{2})} +
  \alpha\frac{\rk(E_{2})}{\rk(E_{1})+
    \rk(E_{2})}.
\end{align*}
If we pick $\alpha$ so that
\begin{equation}\label{eqtn:c=mualpha}
\mu_{\alpha}(\EEE_1,\EEE_2,\Phi)=c
\end{equation}
\noi then a direct calculation shows
\begin{align*}
& \mu(\EEE_1') < c \Longleftrightarrow
\mu_{\alpha}(\EEE'_1,0,0)<\mu_{\alpha}(\EEE_1,\EEE_2,\Phi)\\
& \mu(\EEE_1/\EEE'_1) >c\Longleftrightarrow
\mu_{\alpha}(\EEE'_1,\EEE_2,\Phi)<\mu_{\alpha}(\EEE_1,\EEE_2,\Phi)
\end{align*}
The condition in \ref{thm:fixedE1HK.2}(2) is thus equivalent to the
condition that all subtriples of the form $(\EEE'_1,\EEE_2,\Phi)$ or
$(\EEE'_1,0,\Phi)$ satisfy the
$\alpha$-stability condition for triples with
$\alpha$ defined by (\ref{eqtn:c=mualpha}), i.e. it is
equivalent to the usual triples stability property, but applied only
to subtriples in which $\EEE'_1$ is either the whole  $\EEE_1$ or zero.
\end{remark}
\subsection{Example 3 (Coherent Systems)}\label{subs:CS}
We consider a special case of Example 2, i.e. of the situation in
Section \ref{subs:triples}. Switching notation slightly (to highlight
the different roles played by the two factors) we take
\begin{itemize}
\item $G=\GL(n)\times \GL(k)$,
\item $K=\U(n)\times \U(k)$,
\item $\V=\Hom(\CC^{k},\CC^{n})$, with $\rho:G\longrightarrow \GL(\V)$
given by $\rho(B,C)(T)=B\circ T\circ C^{-1}$.
\end{itemize}
We again take a principal $G$-bundle of the form
$P_G=P_{\GL(n)}\times P_{\GL(k)}$, but we now impose the restriction that
$P_{\GL(k)}$ should be the trivial bundle, i.e.
\begin{equation}
P_{\GL(k)}=X\times \GL(k) \ .
\end{equation}
\noi If, moreover, we fix the trivial holomorphic structure on
$P_{\GL(k)}$, then the remaining data in a pair on
$(P_G,\VV_G=P_G\times_{\rho}\V)$ determines (a) a holomorphic structure on
$P_{\GL(n)}$ and (b) a holomorphic section of $\VV_G$.
Equivalently, replacing the principal bundles with their associated
vector bundles, we see that such pairs correspond to choices for
$(\EEE,\Phi)$, where $\EEE$ is a holomorphic bundle with
$P_{\GL(n)}$ as its frame bundle and $\Phi$ is a holomorphic bundle
map $\Phi:\OOO^{k}\longrightarrow\EEE$. If we fix an identification
\begin{equation}
H^0(\OOO^{k})=\CC^{k}
\end{equation}
then $\Phi$ defines a map from $\CC^{k}$ into $H^0(\EEE)$. We thus
arrive at an interpretation of such principal pairs as {\it coherent
systems}, where a coherent system consists of a holomorphic vector
bundle $\EEE$ together with a linear subspace in $H^0(\EEE)$ (cf.
\cite{KN}, \cite{LeP}). More precisely,
it is given by a vector bundle, $\EEE$,
and a homomorphism $u:\mathbb V\longrightarrow H^0(\EEE)$, where
$\mathbb V$ is a fixed $k$-dimensional vector space.
If we fix an identification $\mathbb V\simeq H^0(\Ok)$, then $u$
defines a map $\Phi:\Ok\rightarrow\EEE$ and the coherent system
$(\EEE,u)$ determines the triple $(\EEE,\Ok, \Phi)$. Conversely, any
such triple determines a coherent system in which $V=H^0(\Ok)$ and
$u$ is the map induced map by $\Phi$. Changing the identification of
$\mathbb V$ with
$H^0(\Ok) $ is equivalent to changing
$u$ by the action of an element in $\GL(\mathbb V)\simeq \GL(k,\CC)$. We thus
get a bijective correspondence between coherent systems modulo the
action of $\GL(k)$ on $V$, and triples $(\EEE,\Ok,\Phi)$ modulo the
action of $\GL(k)$ on $\Ok$. Triples of this sort are thus precisely
the principal pairs described above.
Notice that since $P_{\U(k)}$ is the trivial $\U(k)$-bundle, we can
identify larger subgroups than those which appear in section
\ref{subs:triples}. In particular, as in
\S\ref{sectn:constant} we can define subgroups of
$\GGC(P_G)$ and  $\GGG(P_K)$ respectively by
$$\HHC=\GGG(P_{\GL(n)})\times \GL(k)
\qquad and \qquad\HHH=\GGG(P_{\U(n)})\times \U(k).$$
Moreover, if we take the trivial
holomorphic structure on $P_{\GL(k)}$, then the resulting subspace of
$\XXX_H\subset \chi(P_G,\VV_G)$ is $\HHC$-invariant.
Let $\hlie=\Lie\HHH$. We define the projection
$$\pi_{\hlie}:\Lie\GGG=\Lie(\GGG(P_{\U(n)}))\oplus \Lie(\GGG(P_{\U(k)}))
\to \Lie(\GGG(P_{\U(n)}))\oplus \ulie(k)$$
by integrating the second component over the base manifold $X$, i.e.
\begin{equation}
\pi(g_1,g_2)=\left(g_1,\int_X g_2\right).
\end{equation}
As in the previous examples:
\begin{prop} The subgroups
$\HHH$ and $\HHC$, together with the subspace
$\chi_H$, define a subgroup setting in which our Main Theorem applies.
\end{prop}
In view of the above discussion on the relation between the points in
$\XXX_H$ and coherent systems, we see that the $\HHC$ orbits
in $\XXX_H$ correspond to isomorphism classes of coherent systems.
The result of our Main Theorem, which relates the moment map for the
$\HHH$-action on $\XXX_H\subset \XXX(P_K,\VV_K)$ to stability with respect
to $\HHC$,
can thus be interpreted as a Hitchin--Kobayashi correspondence for
coherent systems.
\begin{prop}\label{thm:CS}
 Let $\HHH$, $\HHC$ and $\XXX_H\subset\XXX(P_G,\VV_G)$ be as above.
Use the standard representations to define an auxiliary
representation
$\rho_a:K\longrightarrow \U(\CC^{n}\oplus\CC^{k})$, and thereby
an inclusion $\Lie\GGG\subset(\Lie\GGG)^*$ (as in
Section \ref{subs:geomsetting}).  Let $c_{\HHH}:=-\sqrt{-1}(c_1I_n,
c_2I_k)$, where $c_1,c_2\in\RR$ are arbitrary real numbers.
Let $((A,0),\Phi)$ be a point in
$\XXX_H$, where $A$ is a connection on $P_{\GL(n)}$ and $0$
denotes the trivial connection on $P_{\GL(k)}$. Fix a global frame
$\{e_i\}_{i=1}^{k}$ for $\OOO^{k}$ and define $\phi_i=\Phi(e_i)$.
Let $S=\mathrm {Span}\{\phi_1,\dots,\phi_k\}$ be the subspace of
$H^0(\EEE)$ spanned by the image of the induced map
$\Phi:H^0(\OOO^k)\longrightarrow H^0(\EEE)$.
\begin{enumerate}
\item When applied to $((A,0),\Phi)$, the
$(\HHH,c_{\HHH})$-vortex equation takes the form
\begin{equation}\label{eqtn:CS}
\begin{split}
\sqrt{-1}\Lambda F_{A}+\Phi\Phi^* = c_1I_n\\
\la\phi_i,\phi_j\ra_{L^2} = -c_2I_{k}
\end{split}
\end{equation}
\noi where the inner product $\la\ ,\ \ra_{L^2}$ is on $H^0(\EEE)$, and
$I_{k}$
denotes the $k\times k$ identity matrix. There are no solutions
unless $c_1$ and $c_2$ satisfy the constraint
\begin{equation}\label{eqn:constraint}
\deg(\EEE)=c_1\rk(\EEE)+c_2k\ .
\end{equation}
\item Assume that $c_1$ and $c_2$ satisfy the constraint
(\ref{eqn:constraint}).
Then a point
$((A,0),\Phi)\in\XXX_H$ is $(\HHH_H,c_{\HHH})$-stable if and only
if for any coherent subsheaf
$\EEE'\subset\EEE$:

\begin{equation}\label{eqn:CSstability}
\frac{\deg(\EEE')}{\rk(\EEE')}+
 \alpha\frac{k'}{\rk(\EEE')} < c_1 \ ,
\end{equation}

 \noi where $k'=\mathrm {dim}(H^0(\EEE')\cap S)$ and $\alpha$ is
 defined by
 \begin{equation}\label{eqn:alphadefn}
 \frac{\deg(\EEE)}{\rk(\EEE)}+
 \alpha\frac{\mathrm {dim}(S)}{\rk(\EEE)} = c_1\ .
 \end{equation}
\end{enumerate}
\end{prop}
\begin{pf} We have to compute $\pi_{\hlie}(\Lambda F_{(A,0)}+\mu(\Phi))$.
If we fix a frame $\{e_i\}$ for $\OOO^{k}$, and define sections of $\EEE$ by
\begin{equation}
\phi_i=\Phi(e_i)
\end{equation}
\noi then with respect to this frame, the endomorphism
$\Phi^*\Phi$ is given by the matrix whose $(ij)$ element is
$\phi_i\overline{\phi_j}$. We thus find that
$$\int_X\sqrt{-1}\Phi^*\Phi=\sqrt{-1}\la\phi_i,\phi_j\ra_{L^2}.$$
At a point $((A,0),\Phi)\in \chi$, where $A$ is a connection on $P_{\GL(n)}$
and $0$ denotes the trivial connection on $P_{\GL(k)}$, we thus get
\begin{eqnarray}
\pi_{\hlie}(\Lambda
F_{(A,0)}+\mu(\Phi))& =& \pi_{\hlie}(\Lambda F_{A}, 0) +
\pi_{\hlie} (-\sqrt{-1}\Phi\Phi^*\ ,\ +\sqrt{-1}\Phi^*\Phi\ )\\
& =& (\Lambda F_{A}-\sqrt{-1}\Phi\Phi^*)\ ,
\ \sqrt{-1}\la\phi_i,\phi_j\ra_{L^2}).
\end{eqnarray}
\noi The equations (\ref{eqtn:CS}) follow from this. The constraint
on $c_1$ and $c_2$ is obtained by integrating the trace of the first
equation in (\ref{eqtn:CS}) and adding the result to the trace of the
second equation.
To prove that the stability conditions one has to follow the same
scheme as in Proposition \ref{prop:fixedE1HK}.
By our choice of auxiliary representation, the vector bundle
$W=P_K\times_{\rho_a}(\CC^n\oplus \CC^k)$ can be identified as $W=\EEE\oplus
\OOO^k$. Any holomorphic filtration induced by $\chi\in\Omega^0(\ad P_K)$ is
of
the form
$$0=W^0\subseteq W^1\subseteq\dots\subseteq W^r=W$$
with $W^j=\EEE^j\oplus\FFF^{k_j}$, where $\EEE^j\subset\EEE$ are
subbundles defining a filtration
$$0=\EEE^0\subseteq \EEE^1\subseteq\dots\subseteq \EEE^r=\EEE$$
and $\FFF^{j}\subset\OOO^k$ are subbundles defining a filtration
$$0=\FFF^{0}\subseteq \FFF^{1}\subseteq\dots\subseteq \FFF^r=\OOO^k\ .$$
Suppose now that $\chi$ is in $\hlie$. Then $\FFF=\OOO^{k_j}$ for
some $k_j\le k$, so that the filtration of $\OOO^k$ is of the form
$$0=\OOO^{k_0}\subseteq \OOO^{k_1}\subseteq\dots\subseteq \OOO^k$$.
Hence, if the eigenvalues of $-\sqrt{-1}\chi$ are
$\alpha_1\le\alpha_2\le\dots\alpha_r$, then
\begin{equation}
\deg\chi=\alpha_r\deg \EEE^1+\sum_{j=1}^{r-1}(\alpha_j-\alpha_{j+1})
\deg\EEE^j.
\end{equation}
Notice however that with $c_{\HHH}:=-\sqrt{-1}(c_1I_n, c_2I_k)$, (and
$\Vol X=1$), we get
\begin{equation}
\int_X\la\chi_{\alpha},c_{\HHH}\ra=
c_1(\alpha_r\rk \EEE+\sum_{j=1}^{r-1}(\alpha_j-\alpha_{j+1})
\rk\EEE_{j})+c_2(\alpha_r k+\sum_{j=1}^{r-1}(\alpha_j-\alpha_{j+1}){k_j})\ .
\end{equation}
\noi The eigenvalues for $-\sqrt{-1}\rho(\chi)$ on $V=\Hom(\OOO^k,\EEE)$
are $\alpha_i-\alpha_j$.
As in the proof of Proposition \ref{prop:fixedE1HK}, we can apply
(SSC) and show that it is enough to consider $\chi\in\hlie$ with at
most two eigenvalues, all of which come from the set $\imag\{-1,0,1\}$. We
may thus assume that the corresponding filtrations are of the form
$$ 0\subseteq \EEE'\oplus\OOO^{k'}\subseteq\EEE\oplus\OOO^k\ ,$$
and the eigenvalues are either $(-\imag,0)$ or $(0,\imag)$. In both cases
the
eigenvalues for $-\sqrt{-1}\rho(\chi)$ on $V=\Hom(\OOO^k,\EEE)$ are
$(-1,0,1)$ and the condition $\Phi\in H^0(V^-)$ is equivalent to
$\Phi(\OOO^{k'})\subset\EEE'$. If the eigenvalues are $(-1,0)$ then
the stability condition $\deg\chi -
\int_X\la\chi_{\alpha},c_{\HHH}\ra>0$ is equivalent to
\begin{equation}\label{eq:case1}
\deg(\EEE')-c_1\rk(\EEE')-c_2k'<0.
\end{equation}
If the eigenvalues are $(0,\imag)$ then the condition $\deg\chi -
\int_X\la\chi_{\alpha},c_{\HHH}\ra>0$ is equivalent to
\begin{equation}\label{eq:case2}
\deg(\EEE')-c_1\rk(\EEE')-c_2k'<\deg(\EEE)-c_1\rk(\EEE)-c_2k.
\end{equation}
However by (\ref{eqn:constraint}), $\deg(\EEE)-c_1\rk(\EEE)-c_2k=0$.
Thus (\ref{eq:case2}) is the same as (\ref{eq:case1}). Defining
$\alpha$ as in (\ref{eqn:alphadefn}), and using (\ref{eqn:constraint}),
it follows immediately that this condition is the same as
(\ref{eqn:constraint}).
\end{pf}

\subsection{Example 4 (Twisted Triples)}
Examples 1 and 2 can be combined to describe objects which might be
called twisted triples; by this we mean a fixed triple, but in which
one of the bundles is only partially fixed. Such objects arise
naturally as a result of dimensional reduction when the reduction is
on a `twisted' orbit space, as in \cite{BGK}. They are also the
simplest example of the twisted quiver bundles considered in
\cite{AG}.  Higgs bundles may be viewed as a special case.
\begin{definition}\label{defn:Ttriple} A twisted triple consists of
two holomorphic bundles $\EEE_1$ and $\EEE_2$, plus a bundle map

 \begin{equation}
 \label{eqtn:Phi1}
 \Phi:\EEE_2\otimes\FFF\longrightarrow \EEE_1
 \end{equation}

 \noi where $\FFF$ is a fixed (`twisting') holomorphic bundle.
\end{definition}
\noi Writing $\tilde{\EEE}_2=\EEE_2\otimes\FFF$, we can think of the
twisted triple as a triple $(\EEE_1,\tilde{\EEE_2},\Phi)$ in which
the second bundle is partially fixed.  Alternatively, we can replace
$\Phi$ by either of the equivalent maps (which, by abuse of notation
we also denote by $\Phi$)
\begin{equation}
\label{eqtn:Phi2}
 \Phi:\EEE_2\longrightarrow \EEE_1\otimes\FFF^*
\end{equation}
\noi or
\begin{equation}
\label{eqtn:Phi3}
 \Phi:\FFF\longrightarrow \EEE_1\otimes\EEE^*_2\ .
\end{equation}
\noi The twisted triple is then be viewed as a triple
$(\EEE_1\otimes\FFF^*,\EEE_2^*,\Phi)$ in which the
first bundle is partially fixed, or as a triple
$(\EEE_1\otimes\EEE_2^*,\FFF,\Phi)$ in which the second bundle is
fixed and the first bundle is a tensor product.
Adopting the first description (i.e. the one given in Definition
\ref{defn:Ttriple}), we get a description in terms of principal
pairs if we take $p=3$ in (\ref{eqtn:Gproduct}) and
(\ref{eqtn:PGproduct}) and consider the case where
\begin{itemize}
\item $G=\GL(n_1)\times \GL(n_2)\times \GL(n_3)$,
\item $K=\U(n_1)\times \U(n_2)\times \U(n_3)$,
\item $\V_1=\CC^{n_1}$ and $\rho_1:\GL(n_1)\longrightarrow \GL(\V_1)$
 is the standard representation, and
\item for $i=2,3$, $\V_i=\CC^{n_i}$ and $\rho_i:\GL(n_i)\longrightarrow
\GL(\V_i)$
 is the dual of the standard representation.
\end{itemize}
We take $ \V=\V_1\otimes \V_2\otimes \V_3$ and let
$\rho=\rho_1\otimes\rho_2\otimes\rho_3$  be the tensor product
representation of $G$ on $\V$.\footnotemark
\footnotetext{This is equivalent to setting
$\V=\Hom(\CC^{n_2}\otimes\CC^{n_3},\CC^{n_1})$ and taking the representation
$\rho(C_1,C_2,C_3)(T)=C_1\cdot T\cdot (C_2\otimes C_3)^{-1}$}
As in Example 1 (in section
\ref{subs:tensor}) $\VV_G=P_G\times_{\rho}\V$ is then a tensor
product of vector bundles; in this case
$\VV_G=\VV_1\otimes\VV^*_2\otimes \VV^*_3=\Hom(\VV_2\otimes\VV_3,\VV_1)$,
where $\VV_i$ denotes the standard vector bundle associated to the
principal $\GL(n_i)$ bundle $P_{\GL(n_i)}$, and $V_i^*$ is the dual of
$V_i$.  A holomorphic structure on
$P_G$ is thus equivalent to a holomorphic structure on
$\VV_G$ such that the resulting holomorphic vector bundle
$\VVV=\VVV_1\otimes\VVV^*_2\otimes \VVV^*_3$ (where $\VVV_i$ denotes
holomorphic bundle obtained by putting a holomorphic structure on
$\VV_i$).
The data in a principal pairs of this type is clearly equivalent to
the defining data in a twisted triple, but without any extra
conditions imposed on the `twisting bundle'. Indeed we can make the
correspondence explicit if we let $\EEE_1$ and $\EEE_2$ be the
standard (holomorphic) vector bundles associated to $P_{\GL(n_1)}$ and
$P_{\GL(n_2)}$ respectively, and let the 'twisting bundle' $\FFF$ be
the standard vector bundle associated to
$P_{\GL(n_3)}$. Depending on whether $\rho_i$ is the standard
representation or its dual, $\VVV_i$ is identified with either
$\EEE_i$ or its dual $\EEE^*_i$, and the section $\Phi$  becomes the
map $\Phi$ in (\ref{eqtn:Phi1}), (\ref{eqtn:Phi2}) or
(\ref{eqtn:Phi3}).
To impose the constraint on $\FFF$, we consider
the normal subgroup of $K$ defined by
\begin{equation}
H=\U(n_1)\times \U(n_2)\times {1} \ .
\end{equation}
\noi Let $\HHH$ be the corresponding subgroup of
$\GGG$, as in Section \ref{sect:normal} and let $\HHC$ be its
complexification. If we fix a connection, say $A_f^0$, on
$P_f=P_{\U(n_3)}$, then we get a subspace
$\XXX_H\subset\XXX(P_K,\VV_K)$ defined by the condition that
the connections on $P_K$ is of the form
\begin{equation}
A=(A_1,A_2, A^0_f)\ .
\end{equation}
\begin{lemma}
The subspace $\XXX_H$ is $\HHH$-invariant. The corresponding complex
subspace of $\XXX(P_G,\VV_G)$ is $\HHC$-invariant and defines the
configuration space of twisted triples. The $\HHC$ orbits correspond
to the isomorphism classes of twisted triples. The set
$(\HHH,\HHC,\XXX_H)$ determines a subgroup setting.
\end{lemma}
\begin{prop}\label{prop:TTriples} Let $\HHH\subset \GGG$and
$\XXX_H\subset\XXX(P_K,\VV_K)$ be as above. Define an auxiliary
representation
$\rho_a:K\longrightarrow \U(\CC^{n_1}\oplus\CC^{n_2}\oplus\CC^{n_3})$ using
the standard representations, and thereby fix an inclusion
$\Lie\GGG\subset(\Lie\GGG)^*$ (as in Section
\ref{subs:geomsetting}). Fix a central element
$c_{\HHH}:=-\sqrt{-1}(c_1I_1,c_2I_2,c_3I_3)$ in the Lie algebra of
$\HHH$ and assume that
$$n_1c_1+n_2c_2=\deg(\EEE_1)+\deg(\EEE_2)\ ,$$
\noi where $\EEE_1$ and $\EEE_2$ are the vector bundles
associated to the principal $\GL(n_1)$
and $\GL(n_2)$ bundles.
(1) Fix a system of local unitary frames for $\VV_3$, say
$\{e^{\alpha}_i\}_{i=1}^{n_3}$ over $U_{\alpha}\subset X$, where
$\{U_{\alpha}\}$ is a suitable open cover of $X$. Let
$\{f^{\alpha}_i\}_{i=1}^{n_3}$ be the dual frame for $\VV_3^*$.
For any $x\in U_{\alpha}$ write
$$\Phi(x)=\sum_{i=1}^{n_3}\Phi^{\alpha}_i(x)\otimes
f^{\alpha}_i(x)$$
\noi where the $\Phi^{\alpha}_i(x)$ are locally defined
sections of $\VV_1\otimes\VV^*_2=\Hom(\VV_2,\VV_1)$. Then when applied to
any
$((A_1, A_2,A_f),\Phi)\in\XXX_H$, the
$(\HHH,c_{\HHH})$-vortex equation takes the form
\begin{eqnarray}
\sqrt{-1}\Lambda F_{A_1}+\sum_{j=1}^{n_f}\Phi_j\Phi_j^*=c_1I_1\\
\sqrt{-1}\Lambda F_{A_2}-\sum_{j=1}^{n_f}\Phi_j^*\Phi_j=c_2I_2.
\end{eqnarray}
\noi The $\Phi_j$ (here viewed as elements in $\Hom(\VV_2,\VV_1)$)
depend on the choice of local frame for
$\FFF=\VVV_3$, but the expressions $\sum_{j=1}^{n_f}\Phi_j\Phi_j^*$
and $\sum_{j=1}^{n_f}\Phi_j^*\Phi_j$ do not.
(2)  Let $(\EEE_1,\EEE_2\otimes\FFF,\Phi)$ be the twisted triple
corresponding to the point $((A_1, A_2,A_f),\Phi)$ in $\XXX_H$. Then
$((A_1, A_2,A_f),\Phi)$ is $(\HHH, c_{\HHH})$-stable if and only if for
every choice of coherent subsheaves $\EEE'_1\subset\EEE_1$ and
$\EEE'_2\subset\EEE_2$ such that $\Phi(\EEE'_2\otimes\FFF)\subset\EEE'_1$,
the
inequality
 $$\frac{\deg(\EEE_1'\oplus\EEE'_2)}{\rk(\EEE_1'\oplus\EEE'_2)}+
 \alpha\frac{\rk(\EEE'_2)}{\rk(\EEE_1'\oplus\EEE'_2)} <
 \frac{\deg(\EEE_1\oplus\EEE_2)}{\rk(\EEE_1\oplus\EEE_2)}+
 \alpha\frac{\rk(\EEE_2)}{\rk(\EEE_1\oplus\EEE_2)} $$

\noi is satisfied. Here $\alpha$ is determined by the identity
$$\frac{\deg(\EEE_1\oplus\EEE_2)}{\rk(\EEE_1\oplus\EEE_2)}+
 \alpha\frac{\rk(\EEE_2)}{\rk(\EEE_1\oplus\EEE_2)}=c_1\ .$$
\end{prop}
\begin{pf} The form of the $(\HHH,c_{\HHH})$ vortex equations follows in the
same was as in Examples 1 and 2. For more details on the precise form
of the terms involving $\Phi$ see the Appendix. The characterization
of the stability condition follows by exactly the same methods as in
the previous examples; it is a very minor modification of the
calculations in the proofs of Proposition \ref{prop:fixedE1HK} and
Proposition \ref{thm:CS}. We thus omit the details.
\end{pf}

\noi Our Main Theorem thus reduces to the Hitchin--Kobayashi correspondence
for
twisted triples, as in \cite{BGK}.
\begin{remark}
\begin{enumerate}
\item If $\FFF$ is the trivial rank $k$ bundle $\Ok$, then we can fix a
global frame for it. The maps
$\Phi_j:\EEE_2\rightarrow \EEE_1$ are then globally defined and the twisted
triple
$(\EEE_1,\EEE_2\otimes\FFF,\Phi)$ can equivalently be described as a
pair of bundles together with $k$ prescribed maps between them.
\item The example of the twisted triples can be generalized in various ways
if we allow $p>3$, i.e. if we admit more than three bundles. Some
such examples are useful for the description of quiver
representations.  The computation of the appropriate moment map is
essentially the same as that in Proposition \ref{prop:TTriples} but
involves more complicated notation. The key is the Moment Map Lemma
(with subsequent generalizations) given in the Appendix. The
description of the stability condition is likewise more complicated
but involves no new ideas.
\end{enumerate}
\end{remark}
\subsection{Example 5 ($\GL(m)$-Higgs bundles)}
In this example we consider the case where
\begin{itemize}
\item $G=\GL(m)\times \GL(n)$ and $K=\U(m)\times \U(n)$,
 with $n=dim_{\bf C}(X)$,
\item $\V_1=\Hom(\CC^{m},\CC^{m})$, and $\rho_1:\GL(m)\longrightarrow
\GL(\V_1)$ is the representation
$$\rho_1(A)B=ABA^{-1}\ ,$$
\item $\V_2=\CC^n$, and $\rho_2:\GL(n)\longrightarrow \GL(\V_2)$ is the
standard representation on $\CC^n$.
\end{itemize}
\noi We also stipulate that
\begin{itemize}
\item $P_{\GL(n)}$ is the frame bundle for the holomorphic cotangent
bundle $(T^*_X)^{1,0}$ (thus, as a smooth bundle, we can identify
$P_{\GL(n)}\times_{\rho_2}\CC^m=(T^*_X)^{1,0}$).
\end{itemize}
As in the previous examples, we take $ \V=\V_1\otimes \V_2$ and let
$\rho=\rho_1\otimes\rho_2$  be the tensor product
representation of $G$ on $\V$.  The data set which determines a
principal pair on
$(P_G,\VV_G=P_G\times_{\rho}\V)$ is thus equivalent to the data in
$(\EEE,\TTT,\Theta)$, where $\EEE$ is a holomorphic bundle with the
topological type of the standard vector bundle associated to
$P_{\GL(m)}$, $\TTT$ is a holomorphic bundle with the topological type of
$(T^*_X)^{1,0}$, and $\Theta$ is a holomorphic section of
$\End(\EEE)\otimes\TTT$.  Suppose now that the holomorphic structure on
$P_{\GL(n)}$ is chosen such that the identification
$P_{\GL(n)}\times_{\rho_2}\CC^n=(T^*_X)^{1,0}$ is an identification
of holomorphic bundles.
{\it With $\TTT$ fixed to be $(T^*_X)^{1,0}$, the remaining data in
principal pairs on
$(P_G,\VV_G)$  determine pairs $(\EEE,\Theta)$, where $\Theta$ is a
holomorphic section of $\End(\EEE)\otimes(T^*_X)^{1,0}$, i.e. the
principal pairs correspond to Higgs bundles on $P_{\GL(m)}$.
Equivalently, if $E$ is the standard vector bundle associated to
$P_{\GL(m)}$, we say that $(\EEE,\Theta)$ defines a Higgs bundle on
$E$.}
To impose the constraint on
$\TTT$, we consider the normal subgroup of $K$ defined by
\begin{equation}
H={\U(m)}\times  1 \ .
\end{equation}
\noi Let $\HHH$ be the corresponding subgroup of
$\GGG$, as in Section \ref{sect:normal} and let $\HHC$ be its
complexification.
If we fix a connection, say $A^0$, on
$P_{\U(n)}$, then we get a subspace
$\XXX_H\subset\XXX(P_K,\VV_K)$ defined by the condition that
the connections on $P_K$ is of the form
\begin{equation}
A=(A,A^0)\ .
\end{equation}
\begin{lemma}
The subspace $\XXX_H$ is $\HHH$-invariant. The corresponding complex
subspace of $\XXX(P_G,\VV_G)$ is $\HHC$-invariant and defines the
configuration space of Higgs bundles on $E$. The $\HHC$ orbits
correspond to the isomorphism classes of Higgs bundles. The set
$(\HHH,\HHC,\XXX_H)$ determines a subgroup setting.
\end{lemma}
\begin{prop}\label{prop:Higgsmmap} Use the
standard representations to define an auxiliary representation
$\rho_a=\rho_{a,m}\oplus\rho_{a,n}:K\longrightarrow
\U(\CC^{m}\oplus\CC^{n})$, and thereby
an inclusion $\Lie\GGG\subset(\Lie\GGG)^*$ (as in
Section \ref{subs:geomsetting}). Identifying $\Lie\HHH$ with
$\Lie\GGG(P_{\U(m)})$, the
$\HHH$-moment map, restricted to points in $\XXX_H$, is given by
\begin{equation}
\label{eqtn:Higgsmmap}
\mu_H(A,A^0,\Theta)=\Lambda F_A
 -
\sqrt{-1}\Lambda [\Theta, \overline{\Theta}^t].
\end{equation}
\end{prop}
Here $\overline{\Theta}^t$ is defined as follows. Given any $x\in X$,
we can fix a local unitary frame, say $\{e_i\}_{i=1}^n$, for
$(T^*_X)^{1,0}$. We can then write $\Theta=\sum_{i=1}^n\Theta_i\otimes
e_i$, where the $\Theta_i$ are locally defined sections of $\End(E)$.
Then
\begin{equation}
\overline{\Theta}^t =
\sum_{i=1}^n\Theta^*_i\otimes\overline{e_i}
\end{equation}
\noi where the adjoints $\Theta^*_i$ are with respect to the metric
on $E$.

\begin{pf}{\it (of Proposition \ref {prop:Higgsmmap}})
The first term, i.e. $\Lambda F_A$, is obtained in the usual way as a
result of the projection from the Lie algebra of the full gauge group
onto the Lie algebra of $\HHH$. The second term can be understood as
follows. The moment map for the
$\HHH$-action on $\Omega(\End(\EEE)\otimes(T^*_X)^{1,0})$ is,
by Proposition \ref{prop:muVK} in the Appendix,
\begin{equation}
\mu(\Theta)=\sum_{i=1}^n\mu_1(\Theta_i)
\end{equation}
\noi where $\mu_1$ is the moment map for the action of $\U(m)$ on
$\Omega(\End(\EEE)\otimes(T^*_X)^{1,0})$. But, since the moment map
for the conjugation action of $\U(m)$ on $\Hom(\CC^m,\CC^m)$ is
\begin{equation}
\mu_1(A)=-\sqrt{-1}[A,A^*]\ ,
\end{equation}
\noi it follows that
\begin{equation}
\mu(\Theta)=-\sqrt{-1}\sum_{i=1}^n [\Theta_i, \Theta^*_i]\ .
\end{equation}
\noi On the other hand, remembering that $\Theta$ behaves like a 1-form, we
get
\begin{eqnarray}
\Lambda [\Theta, \overline{\Theta}^t]& =&
\sum_{j=1}^n\sum_{i=1}^n \Theta_i\Theta^*_j\Lambda
e_i\wedge\overline{e_j}+\Theta^*_j\Theta_i\Lambda
\overline{e_j}\wedge e_i \\
& =& \sum_{j=1}^n \Theta_i\Theta^*_i-\Theta^*_i\Theta_i\\
& =& \sum_{j=1}^n [\Theta_i,\Theta^*_i]\ ,
\end{eqnarray}
\noi where we have used the fact that
$\Lambda e_i\wedge\overline{e_j}=-\Lambda \overline{e_i}\wedge e_j=0$
if $i\ne j$, and also that
$\Lambda e_i\wedge\overline{e_i}=
-\Lambda \overline{e_i}\wedge e_i=1$.
\end{pf}

\noi We thus get, as an immediate corollary:
\begin{corollary}\label{cor:Higgsvortex}
Given a central element $c_{\HHH}=-\sqrt{-1}(0,c_mI_m)$ in
$\Lie\HHH$,
\begin{enumerate}
\item a point in
$\XXX_H$ thus satisfies the
$(\HHH, c_{\HHH})$-vortex equations if and only if it satisfies the
conditions
\begin{equation}
\label{eqn:Higgseqtn}
\sqrt{-1}\Lambda F_A+\Lambda [\Theta, \overline{\Theta}^t]=c_mI_m\ ,
\end{equation}
\item there are no solutions unless $c_m=\mu(\EEE)$.
\end{enumerate}
\end{corollary}
\begin{proof} (1) follows immediately from proposition
\ref{prop:Higgsmmap}. (2) follows from (1) by integrating the trace of
the (\ref{eqn:Higgseqtn}) and observing that  the trace of $[\Theta,
\overline{\Theta}^t]$ is zero.
\end{proof}
\noi Furthermore,
\begin{prop}\label{prop:Higgstability} A point in $\XXX_H$ is
$(\HHH, c_{\HHH})$-stable if and only if the
corresponding Higgs bundle $(\EEE,\Theta)$ satisfies the condition
$$\mu(\EEE')<\mu(\EEE)$$
\noi for all $\Theta$-invariant coherent subsheaves $\EEE'\subset\EEE$.
\end{prop}
\begin{pf}
We proceed as in the proof of Proposition \ref{prop:fixedE1HK}. With
auxiliary representation $\rho_a$ as above in Proposition
\ref{prop:Higgsmmap}, we consider elements
$\chi_{\alpha}=(\chi_{1,\alpha},0)\in\hlie$ where $\chi_{1,\alpha}$
has eigenvalues
$-\sqrt{-1}\{\alpha_1,\dots,\alpha_r\}$ and corresponding
$A$-holomorphic filtrations
$$0\subset\EEE_{1}\subset\dots\subset\EEE_{r}$$
of $\EEE=P_{\U(m)}\times_{\rho_{a,m}}\CC^{m}$. The eigenvalues for
$\rho(\chi)$ on $V=\End(\EEE)\otimes(T^*_X)^{1,0})$ are then the
differences $-\sqrt{-1}(\alpha_i-\alpha_j)$. The negative subbundle
$V^-$ is thus determined by the condition
$\alpha_i\le\alpha_j$. It follows that the condition
$\Theta\in H^0(V^-)$ is equivalent to the condition
$\Theta(\EEE_j)\subset\EEE_j\otimes (T^*_X)^{1,0}$ for $1\le j\le r$.
We now identify filtrations $\chi(f_i)$ and $\chi(g_i)$ and define
the set $S\subset \hlie$ in exactly the same way as in proof of
Proposition \ref{prop:fixedE1HK}. The verification that (SSC1)
applies follows precisely as before. The verification of (SSC2')
follows from the above characterization of the condition
$\Theta\in H^0(V^-)$.
Suppose now that $\sigma\in S$ defines a filtration
$0\subset\EEE'\subset\EEE$. Whether the eigenvalues of $\sqrt{-1}\sigma$ are
$\alpha_1=-1$ and $\alpha_2=0$,  or $\alpha_1=0$ and $\alpha_2=-1$,
the eigenvalues of $\sqrt{-1}\rho(\sigma)$ on
$V=\End(\EEE)\otimes(T^*_X)^{1,0})$ are $(-1, 0, 1)$. The condition
$\Theta\in H^0(V^-)$ is equivalent to the condition
$\Theta(\EEE')\subset\EEE'\otimes (T^*_X)^{1,0}$.
If the eigenvalues of $\sqrt{-1}\sigma$ are $\alpha_1=-1$ and
$\alpha_2=0$, then $\deg(\sigma)>0$ is equivalent to the condition
$$\mu(\EEE')<c_m\ ,$$
while if the eigenvalues are $\alpha_1=0$ and
$\alpha_2=1$ then the condition is equivalent to
$$\mu(\EEE/\EEE')>c_m\ .$$
However by Corollary \ref{cor:Higgsvortex} (2) we may assume
$c_m=\mu(\EEE)$, and hence both conditions are equivalent to
$$\mu(\EEE')<\mu(\EEE)\ .$$
This completes the proof.
\end{pf}

Combining Corollary \ref{cor:Higgsvortex} and Proposition
\ref{prop:Higgstability} our Main Theorem thus becomes the usual
Hitchin--Kobayashi correspondence for Higgs bundles. We remark that
Higgs bundles with more general structure groups (cf. \cite{H}) can
be treated in a similar manner, but we leave the details to the
reader.
\section{Appendix: Moment Map Lemma}
For $1\le i\le p$ let $\V_i$ be a complex vector space of dimension
$n_i$. Let $\la\ ,\ \ra_i$ be a hermitian inner product on $\V_i$ and
let $\omega_i$ be the corresponding Kaehler form. Thus
\begin{equation}
\omega_i(x,y)=\frac{1}{2\imag}(\la x,y\ra_i-\la y,x\ra_i)\ .
\end{equation}
Let $\V$ be the tensor product
$\V=\V_1\otimes \V_2\otimes\dots\otimes \V_p$, and let $\la\ ,\ \ra$ be the
hermitian inner product determined by the inner products on the
$\V_i$. Thus
\begin{equation}
\la x_1\otimes x_2\otimes\dots\otimes x_p,
y_1\otimes y_2\otimes\dots\otimes y_p\ra = \Pi_{i=1}^p\la x_i,y_i\ra_i\ .
\end{equation}
Let $\Omega$ be the corresponding Kaehler form. In this Appendix we
compute moment maps for some Hamiltonian actions on $(\V,\Omega)$.
\begin{prop}\label{prop:A1}
Let $\U_i$ be the group of unitary transformations on
$(\V_i,\la\ ,\ \ra_i)$, and use the inner
product $\la\ ,\ \ra$ to identify the
Lie algebra of $\U_i$ with its dual. Then $\U_i$ acts symplectically on
$\V_i$, with moment map
$$\mu_i:\V_i\longrightarrow \Lie(\U_i)=\mathfrak{u}_i$$
\noi given by
\begin{equation}
\label{eqtn:ApndxbasicMmap}
\mu^{V_i}_i(x)=-\sqrt{-1}x\overline{x}^T
\end{equation}
\noi or equivalently
\begin{equation}
\label{eqtn:ApndxbasicMmap2}
\mu^{V_i}_i(x)=-\sqrt{-1}x\otimes x^*.
\end{equation}
\end{prop}
(In (\ref{eqtn:ApndxbasicMmap}) we regard a vector $x\in \V_i$ as a
$n_i\times 1$ column whose entries are the components of $x$ with
respect to a unitary frame for
$\V_i$. Then  the transpose $x^T$ is a $1\times n_i$ row and
$x\overline{x}^T$ is a $n_i\times n_i$ matrix, i.e. an endomorphism of
$\V_i$.
In (\ref{eqtn:ApndxbasicMmap2}) the vector $x^*$ is the element
corresponding to $x$ under the duality $\V_i\cong \V^*_i$ determined
by the metric. Thus $x\otimes x^*$ is in $\V_i\otimes \V^*_i$, which
we can identify with the endomorphisms of $\V_i$.)
If any other group, $K_i$ acts symplectically on $\V_i$ via a
faithful representation
\begin{equation}
\rho_i:K_i\longmapsto \U_i
\end{equation}
then we get a moment map for the $K_i$ action, denoted by
\begin{equation}
\mu^{V_i}_{K_i}:V_i\longmapsto \Lie(K_i)^*\ .
\end{equation}
\begin{prop}
Fix an inner product on
$\Lie(K_i)$ such that
$\rho_{i*}:\Lie(K_i)\longmapsto \mathfrak{u}_i$ is an isometry, and
use this to identify $\Lie(K_i)\cong(\Lie(K_i))^*$. Then we get
\begin{equation}
\mu^{\V_i}_{K_i}=\pi_{K_i}\circ\mu^{\V_i}_i\ ,
\end{equation}
\noi where $\pi_{K_i}:\mathfrak{u}_i\longmapsto \rho_{i*}(\Lie(K_i))$ is
orthogonal projection onto the linear subspace.
\end{prop}
Each group $\U_i$ acts symplectically on $(\V,\Omega)$ via the action
\begin{equation}
A_i(x_1\otimes x_2\otimes\dots\otimes x_p)=x_1\otimes x_2\otimes\dots
A_ix_i\otimes\dots\otimes x_p \ .
\end{equation}
\noi More generally, each group $K_i$ acts symplectically on
$(\V,\Omega)$ via
\begin{equation}
k(x_1\otimes x_2\otimes\dots\otimes x_p)=x_1\otimes x_2\otimes\dots
\rho_i(k)x_i\otimes\dots\otimes x_p \ .
\end{equation}
To describe the moment maps for these actions, it is convenient to
fix unitary bases $\{e^{(i)}_{j}\}_{j=1}^{n_i}$ for each $\V_i$ and
write $X\in \V$ as
$$X=\sum_{i_1,i_2,\dots,i_p}X_{i_1i_2\dots i_p}e^{(1)}_{i_1}\otimes
e^{(2)}_{i_2}\otimes\dots\otimes e^{(p)}_{i_p} \ .$$
\noi For any $1\le i \le p$ we can think of $\V$ as a tensor product
$\widehat{\V_i}\otimes \V_i$, where $\widehat{\V_i}$ is the tensor product
of all the $\V_1,\dots,\V_p$ except for $\V_i$, and write
$X$ as
\begin{equation}
X=\sum_{j=1}^{n_i}X_j\otimes e^{(i)}_{j}\ ,
\end{equation}
\noi where $X_j$ are vectors in $\widehat{\V_i}$.
\begin{prop}\label{prop:muVK}
The moment map
$\mu^{V}_{K_i}:\V\longrightarrow \Lie(K_i)$ is given by
\begin{equation}
\label{eqtn:mmapdual}
\mu^{\V}_{K_i}(X)=\sum_{j=1}^{n_i}\mu^{\V_i}_{K_i}(X_{j})\ .
\end{equation}
\noi In particular, the moment map for the action of $\U_i$ on $\V$ is
\begin{eqnarray}\label{eqtn:mmap}
\mu^{\V}_{i}(X)& =& -\sqrt{-1}\sum_{j=1}^{n_i}X_j\overline{X_j}^t\\
& =& -\sqrt{-1}\sum_j\sum_{\hat{j}} x_{i_1i_2\dots j\dots i_p}
\overline{x}_{i_1i_2\dots j\dots i_p} e^{(i)}_j\otimes
(e^{(i)}_j)^*\ ,
\end{eqnarray}
\noi where the sum in $\sum_{\hat{j}}$ is over all the indices
{\it except} the j'th one.
\end{prop}
\begin{remark} Notice that, while the definition of the $X_{j}$ depends on
the choices of unitary bases, the combination
$\sum_{j=1}^{n_i}\mu^{\V_i}_{K_i}(X_{j})$ does not.
\end{remark}
\noi We consider a few special cases, which come up in interesting
examples.
\subsection{$p=1$}
In this case we recover the basic moment map given in
(\ref{eqtn:ApndxbasicMmap}).
\subsection{$p=2$ (triples)}
Suppose that $\V=\V_1\otimes \V^*_2$. We can identify
\begin{equation}
\V_1\otimes \V^*_2 =\Hom(\V_2,\V_1)
\end{equation}

and interpret any element
$\Phi=\sum_{i,j}\Phi_{ij}e^{(1)}_{i}\otimes (e^{(2)}_j)^*$ as a map
\begin{equation}
\Phi: \V_2\longrightarrow \V_1\ .
\end{equation}
\noi Indeed the map is given by
\begin{equation}
(e^{(2)}_j)\mapsto\sum_i\Phi_{ij}e^{(1)}_{i}\ .
\end{equation}
\begin{lemma}
Under this identification the $\U_1$ moment map given by
(\ref{eqtn:mmap})becomes
\begin{equation}
\label{eqtn:mmap1}
\mu_1(\Phi) = -\sqrt{-1}\Phi\Phi^*.
\end{equation}
\end{lemma}

We can also compute the moment map for $\U_2$. Switching the roles of
$\V_1$ and $\V^*_2$, and using the dual action on $\V^*_2$, we get
\begin{equation}
\label{eqtn:mmap2}
\mu_2(\Phi) = \sqrt{-1}\overline{\Phi}^T\Phi = \sqrt{-1}\Phi^*\Phi.
\end{equation}
The transpose comes from writing
\begin{equation}
\Phi=\sum_{i,j}\Phi_{ij}e^{(1)}_{i}\otimes (e^{(2)}_j)^*  =
\sum_{i,j}\Phi_{ji}(e^{(2)}_{i})^*\otimes e^{(1)}_j
\end{equation}
while the change of sign and conjugation come from the dual group
action.
Notice that the  moment maps (\ref{eqtn:mmap1}) and
(\ref{eqtn:mmap2}) are none other than the projections onto
$\mathfrak{u}_1$ and $\mathfrak{u}_2$ respectively of the moment map for the $\U_1\times \U_2$
action on $\Hom(\V_2,\V_1)$.
\subsection{$p=3$ (twisted triples)}

Suppose that $\V=\V_1\otimes \WW\otimes \V^*_2$. Under the
identification
\begin{equation}
\V_1\otimes \WW\otimes \V^*_2 =\Hom(\V_2,\V_1\otimes \WW)
\end{equation}
\noi we can interpret vectors in $\V=\V_1\otimes \WW\otimes \V^*_2$
as maps
\begin{equation}
\Phi: \V_2\longrightarrow \V_1\otimes \WW\ .
\end{equation}
\noi Indeed, writing
\begin{equation}
\Phi=\sum_{i,j,k}\Phi_{ijk}e^{(1)}_{i}\otimes f_j\otimes (e^{(2)}_k)^*\ ,
\end{equation}
\noi where $\{f_j\}$\ is a frame for $\WW$, the corresponding
map is given by
\begin{equation}
(e^{(2)}_k)\mapsto\sum_{i,j}\Phi_{ijk}e^{(1)}_{i}\otimes f_j\ .
\end{equation}
Defining $\phi_j\in \V_1\otimes \V^*_2$ by
\begin{equation}
\phi_j=\sum_{i,k}\Phi_{ijk}e^{(1)}_{i}\otimes (e^{(2)}_k)^*
\end{equation}
\noi we can write
\begin{equation}
\Phi=\sum_{j}\phi_j\otimes f_j\ .
\end{equation}
\begin{lemma}
Under these identifications the $\U_1$ and $\U_2$ moment map given by
(\ref{eqtn:mmap}) and (\ref{eqtn:mmapdual}) become
\begin{eqnarray}
\mu_1(\Phi) = -\sqrt{-1}\sum_j\phi_j\phi^*_j\\
\mu_2(\Phi) = \sqrt{-1}\sum_j\phi^*_j\phi_j.
\end{eqnarray}
\end{lemma}
\noi Notice that, while the definition of the $\phi_j$ depend on
the choices of unitary bases, the endomorphisms
$\sum_j\phi_j\phi^*_j$ and  $\sum_j\phi^*_j\phi_j$ are invariantly
defined. These quantities correspond exactly to the terms which arise
in the coupled twisted vortex equations for twisted triples (cf.
\cite{BGK}).


\end{document}